\documentclass[12pt]{amsart}
\begin{document}
\theoremstyle{plain}
\newtheorem{thm}{Theorem}[section]
\newtheorem*{thm*}{Theorem}
\newtheorem{prop}[thm]{Proposition}
\newtheorem*{prop*}{Proposition}
\newtheorem{lemma}[thm]{Lemma}
\newtheorem{cor}[thm]{Corollary}
\newtheorem*{conj*}{Conjecture}
\newtheorem*{cor*}{Corollary}
\newtheorem{defn}[thm]{Definition}
\theoremstyle{definition}
\newtheorem*{defn*}{Definition}
\newtheorem{rems}[thm]{Remarks}
\newtheorem*{rems*}{Remarks}
\newtheorem*{proof*}{Proof}
\newtheorem*{not*}{Notation}
\newcommand{\npartial}{\slash\!\!\!\partial}
\newcommand{\Heis}{\operatorname{Heis}}
\newcommand{\Solv}{\operatorname{Solv}}
\newcommand{\Spin}{\operatorname{Spin}}
\newcommand{\SO}{\operatorname{SO}}
\newcommand{\ind}{\operatorname{ind}}
\newcommand{\Index}{\operatorname{index}}
\newcommand{\ch}{\operatorname{ch}}
\newcommand{\rank}{\operatorname{rank}}
\newcommand{\abs}[1]{\lvert#1\rvert}
 \newcommand{\A}{{\mathcal A}}
        \newcommand{\D}{{\mathcal D}}\newcommand{\HH}{{\mathcal H}}
        \newcommand{\LL}{{\mathcal L}}
        \newcommand{\B}{{\mathcal B}}
        \newcommand{\K}{{\mathcal K}}
\newcommand{\oo}{{\mathcal O}}
         \newcommand{\PP}{{\mathcal P}}
        \newcommand{\s}{\sigma}
        \newcommand{\coker}{{\mbox coker}}
        \newcommand{\p}{\partial}
        \newcommand{\dd}{|\D|}
        \newcommand{\n}{\parallel}  
\newcommand{\bma}{\left(\begin{array}{cc}}
\newcommand{\ema}{\end{array}\right)}
\newcommand{\bca}{\left(\begin{array}{c}}
\newcommand{\eca}{\end{array}\right)}

\newcommand{\sr}{\stackrel}
\newcommand{\da}{\downarrow}
\newcommand{\tD}{\tilde{\D}}

        \newcommand{\R}{\mathbf R}
        \newcommand{\C}{\mathbf C}
        \newcommand{\h}{\mathbf H}
\newcommand{\Z}{\mathbf Z}
\newcommand{\N}{\mathbf N}
\newcommand{\tto}{\longrightarrow}
\newcommand{\ben}{\begin{displaymath}}
        \newcommand{\een}{\end{displaymath}}
\newcommand{\be}{\begin{equation}}
\newcommand{\ee}{\end{equation}}

        \newcommand{\bean}{\begin{eqnarray*}}
        \newcommand{\eean}{\end{eqnarray*}}
\newcommand{\nno}{\nonumber\\}
\newcommand{\bea}{\begin{eqnarray}}
\newcommand{\eea}{\end{eqnarray}}
\newcommand{\supp}[1]{\operatorname{#1}}
\newcommand{\norm}[1]{\parallel\, #1\, \parallel}
\newcommand{\ip}[2]{\langle #1,#2\rangle}
\setlength{\parskip}{.3cm}
\newcommand{\nc}{\newcommand}
\nc{\nt}{\newtheorem}
\nc{\gf}[2]{\genfrac{}{}{0pt}{}{#1}{#2}}
\nc{\mb}[1]{{\mbox{$ #1 $}}}
\nc{\real}{{\mathbb R}}
\nc{\comp}{{\mathbb C}}
\nc{\ints}{{\mathbb Z}}
\nc{\Ltoo}{\mb{L^2({\mathbf H})}}
\nc{\rtoo}{\mb{{\mathbf R}^2}}
\nc{\slr}{{\mathbf {SL}}(2,\real)}
\nc{\slz}{{\mathbf {SL}}(2,\ints)}
\nc{\su}{{\mathbf {SU}}(1,1)}
\nc{\so}{{\mathbf {SO}}}
\nc{\hyp}{{\mathbb H}}
\nc{\disc}{{\mathbf D}}
\nc{\torus}{{\mathbb T}}
\newcommand{\tk}{\widetilde{K}}
\newcommand{\boe}{{\bf e}}\newcommand{\bt}{{\bf t}}
\newcommand{\vth}{\vartheta}
\newcommand{\CGh}{\widetilde{\CG}}
\newcommand{\db}{\overline{\partial}}
\newcommand{\tE}{\widetilde{E}}
\newcommand{\tr}{\mbox{tr}}
\newcommand{\ta}{\widetilde{\alpha}}
\newcommand{\tb}{\widetilde{\beta}}
\newcommand{\txi}{\widetilde{\xi}}
\newcommand{\hV}{\hat{V}}
\newcommand{\IC}{\mathbf{C}}
\newcommand{\IZ}{\mathbf{Z}}
\newcommand{\IP}{\mathbf{P}}
\newcommand{\IR}{\mathbf{R}}
\newcommand{\IH}{\mathbf{H}}
\newcommand{\IG}{\mathbf{G}}
\newcommand{\CC}{{\mathcal C}}
\newcommand{\CD}{{\mathcal D}}
\newcommand{\CS}{{\mathcal S}}
\newcommand{\CG}{{\mathcal G}}
\newcommand{\CL}{{\mathcal L}}
\newcommand{\CO}{{\mathcal O}}
\nc{\ca}{{\mathcal A}}
\nc{\cag}{{{\mathcal A}^\Gamma}}
\nc{\cg}{{\mathcal G}}
\nc{\chh}{{\mathcal H}}
\nc{\ck}{{\mathcal B}}
\nc{\cl}{{\mathcal L}}
\nc{\cm}{{\mathcal M}}
\nc{\cn}{{\mathcal N}}
\nc{\NN}{{\mathcal N}}
\nc{\cs}{{\mathcal S}}
\nc{\cz}{{\mathcal Z}}
\nc{\sInd}{\sigma{\rm -Ind}}
\newcommand{\la}{\langle}
\newcommand{\ra}{\rangle}

\begin{center}
 \title{The
local index formula in semifinite von Neumann algebras II: The even case}
 \vspace{.5 in}

\author{}

\maketitle

{\bf Alan L. Carey}\\
Mathematical Sciences Institute\\
Australian National University\\
Canberra, ACT. 0200, AUSTRALIA\\
e-mail: acarey@maths.anu.edu.au\\ 
\vspace{.2 in}

{\bf John Phillips}\\
Department of Mathematics and Statistics\\
University of Victoria\\Victoria, B.C. V8W 3P4, CANADA\footnote{Address for
correspondence}\\
e-mail: phillips@math.uvic.ca\\ 
\vspace{.2 in}

{\bf Adam Rennie}\\
School of Mathematical and Physical Sciences\\
University of Newcastle\\
Callaghan, NSW, 2308 AUSTRALIA\\
e-mail: adam.rennie@newcastle.edu.au\\

\vspace{.2 in}

{\bf Fyodor A. Sukochev}\\
School of Informatics and Engineering\\
Flinders University\\
Bedford Park S.A 5042 AUSTRALIA\\
e-mail: sukochev@infoeng.flinders.edu.au\\ 

\vspace{.25 in}

All authors were supported by grants from ARC (Australia) and
NSERC (Canada), in addition the third named author acknowledges a 
University of Newcastle early career researcher grant.

\end{center}

\newpage
\centerline{{\bf Abstract}}
We generalise the even local index formula of Connes and Moscovici to the
case of spectral triples for a $*$-subalgebra $\A$
of a general semifinite von Neumann algebra. 
 The proof is a variant of that for the odd case which appears in Part I.
To allow for algebras with a non-trivial centre we have to establish 
a theory of unbounded 
Fredholm operators in a general semifinite von Neumann algebra and 
in particular prove a generalised McKean-Singer formula.  
\footnote{
von Neumann algebra, Fredholm module, cyclic cohomology, chern character,
McKean-Singer formula.}

\newpage
\allowdisplaybreaks
\section{Introduction}
There have been two new proofs of
the local index theorem in noncommutative geometry of Connes and 
Moscovici \cite{CM}),  by Higson \cite{hig} and, 
for the odd case, by the present authors 
in  part I of this two part
series of papers \cite{CPRS2}.
The novelty in \cite {CPRS2} is consideration of
spectral triples ``inside'' a general 
semifinite von Neumann algebra and 
in the introduction of a new odd cocycle
(in the $(b,B)$ bicomplex of cyclic cohomology)
which provides a substitute in our approach
for the JLO cocycle \cite{Co4} used in \cite{CM}.
 Our new cocycle is reminiscent of, but distinct from,
Higson's `improper cocycle' \cite{hig}.
In subsequent work \cite{CPRS4}, we 
will relate these two cocycles showing how to
obtain a renormalised version of Higson's cocycle from our resolvent cocycle.

The present paper is concerned with two
primary results, the even semifinite local index formula
proved via the even resolvent cocycle and 
a prerequisite, a general theory of Fredholm operators in von Neumann 
algebras which may have non-trivial centre. This extension is essential to 
encompass examples such as arise in the $L^2$ index theorem of Atiyah.
(Other applications are referenced in Part I.)

For a finitely summable even spectral triple with spectral dimension $q$
(see \cite{CPRS2} for the latter terminology) we use 
the even resolvent cocycle to obtain an expression for the
index. The even resolvent cocycle 
is a $(b,B)$ cocycle with values in functions 
defined and holomorphic in a certain half-plane modulo those functions 
holomorphic in a larger half-plane containing the critical point $r=(1-q)/2$.
By taking residues at the critical point as in   \cite{CPRS2} we
prove the even case of a local index formula for smooth finitely summable 
semifinite spectral triples. Thus
as in \cite{CPRS2} we need the property of
 `isolated spectral dimension' to
analytically continue our
 resolvent cocycle term-by-term 
to a deleted neighbourhood of $r=(1-q)/2$. This
 then defines a generalisation of the Connes-Moscovici even residue cocycle
in the finite $(b,B)$ bicomplex.

There remains one gap in our treatment in that we do not prove that the 
residue cocycle represents the Chern character of our semifinite spectral 
triple. This gap will be filled in a subsequent paper, \cite{CPRS4}, as the 
proof is not short.

Our exposition is organised as follows. We assume all of the notation
of the first part \cite{CPRS2} but include additional
 preliminary material, notation
and definitions needed for this paper in Section 2.
Our main theorem starts from a version of the McKean-Singer formula
for the index. However, we
found that Fredholm theory in semifinite von Neumann 
algebras with a non-trivial 
centre did not exist in a form that was suitable for this purpose. In 
particular, the case of an operator which is Fredholm from the range of one
projection to the range of another projection (which is the case of the
McKean-Singer formula) had not been touched in this setting, and is rather 
subtle.
Thus Section 3 establishes such a theory. We note that in this paper we fix a 
faithful normal semifinite trace $\tau$ on our algebra once and for all. 
Thus strictly speaking we deal always with $\tau$-Fredholm theory, and
do not give a full treatment involving centre-valued traces
and related machinery. 
Those expert in all these matters can move straight to Section 4 where we
 state our main theorem, the local
index theorem for even semi-finite spectral triples. 

The main theorem has three parts. The first expresses an index pairing as 
the {\em residue} of the  pairing between the resolvent cocycle and the Chern 
character of a projection. This residue exists with no assumptions concerning 
analytic continuations. The second statement is similar to the first, but the 
index is expressed as the residue of a sum of zeta functions. The third part 
finally assumes that we can analytically continue the individual zeta 
functions, so that we express the index pairing as a sum of residues of zeta 
functions. These residues assemble to form a $(b,B)$ cocycle, called the 
residue cocycle.

The proof has a number of important differences
with that of the odd case
and these are highlighted in subsection 5.1 
where we establish an analytic formula for the even index which 
is the starting point for our proof. The rest of Section 5 contains the 
computations needed to prove the main theorem.  By  Subsection 5.6 we
have enough to prove part $2)$ of the main theorem and  
the index formula of part $3)$. 
To prove part $1)$ and the cohomological part of $3)$, we introduce the 
resolvent cocycle for the even case in Section 6.

We conclude this introduction with some general comments on the 
existing proofs 
of the Local Index theorem which may help put our results in context. 
Connes and 
Moscovici begin with a representative of the Chern character (the JLO cocycle) 
and deform it to obtain the unrenormalised residue cocycle. It is automatically 
a representative of the Chern character, and so an index cocycle. While this 
cocycle can be renormalised, it is unclear to us whether a procedure exists to 
modify the JLO cocycle so that it yields the renormalised version automatically.

Higson {\em writes down} a function valued cocycle, proves that it is an index 
cocycle and then proves it is in the class of the Chern character \cite{hig}. 
The unrenormalised local index theorem follows from Higson's cocycle and the 
pseudodifferential calculus. We show in \cite{CPRS4} that there is a simple 
modification of Higson's cocycle which leads directly to the renormalised 
residue cocycle. 
In this paper, as in \cite{CPRS2}, we begin with an analytic formula 
for the index 
pairing and apply perturbation theory and the pseudodifferential calculus to 
obtain the renormalised residue cocycle directly. As part of this process we 
also obtain a function valued (almost) cocycle similar to Higson's, but with 
superior holomorphy properties. Our cocycles are automatically index cocycles, 
and so we need only show that they are in the class of the Chern character. 
This will be shown in \cite{CPRS4}, closely following Higson's methods.

\section{Definitions and Background}

We adopt the notational conventions of \cite{CPRS2}. Thus ${\mathcal N}$
is semifinite von Neuman algebra acting on a Hilbert space $\HH$ and
$\tau$ is a faithful normal semifinite trace
on $\mathcal N.$
 An even  semifinite
spectral triple $(\A,\HH,\D)$ is given by a 
$*$-algebra $\A\subset \cn$, 
a densely defined unbounded operator $\D$ affiliated with $\mathcal N$ 
on $\HH$ and
in addition to the properties of definition 2.1
of \cite{CPRS2}, has a grading $\gamma\in\cn$ such that 
$\gamma^*=\gamma$, $\gamma^2=1$,  $a\gamma=\gamma a$ for all $a\in\A$ and 
$\D\gamma+\gamma\D=0$.
As in \cite{CPRS2} we deal only with unital algebras
${\mathcal A}$ where the identity of $\mathcal A$ is that of $\cn$.
We write $P=(1+\gamma)/2$ and $\D^+=(1-P)\D P= P^{\perp}\D P$.
The operator $\D^+:\HH^+=P(\HH)\to\HH^-=P^{\perp}(\HH)$ is, as we shall see,
an unbounded Breuer-Fredholm operator.

The numerical index discussed here is the result of a pairing between an even
 $K$-theory
class represented by a projection $p$, 
and an even $K$-homology class represented
by
$(\A,\HH,\D)$, \cite[Chapter III,IV]{Co4}. This point of view also makes sense
in the
general semifinite setting after suitably interpreting $K$-homology 
classes, \cite{CPRS1,CP2}.
The pairing of $(b,B)$ cocycles with $K$-homology classes
is written in the even case as
\be \langle [p],[(\A,\HH,\D)]\rangle=
\langle [Ch_*(p)],[Ch^*(\A,\HH,\D)]\rangle,\label{Indpair}\ee
where $[p]\in K_0(\A)$ is a $K$-theory class with representative $p$ and
$[(\A,\HH,\D)]$ is the $K$-homology class of the even spectral triple 
$(\A,\HH,\D)$.
On the right hand side, $Ch_*(p)$ is the Chern character of $p$, and
$[Ch_*(p)]$
its periodic cyclic
homology class. Similarly $[Ch^*(\A,\HH,\D)]$ is the periodic cyclic
cohomology class of the Chern
character of $(\A,\HH,\D)$. {\em The analogue of Equation (\ref{Indpair}), for
a suitable cocycle associated to $(\A,\HH,\D)$, in the general semifinite case
is part of our main result.}

We refer to \cite{Co4,Lo,CPRS2} for the definition of the $(b,B)$ bicomplex.
The $(b,B)$ Chern 
character of a 
projection in an algebra $\A$ is an even $(b,B)$ cycle with $2m$-th term$, 
m\geq 1$, given by
$$ Ch_{2m}(p)=(-1)^m\frac{(2m)!}{2(m!)}(2p-1)\otimes p^{\otimes 2m}.$$
For $m=0$ the definition is $Ch_0(p)=p$. 

\section{Fredholm Theory in Semifinite von Neumann Algebras}

We need to generalise the real-valued Fredholm index theory outlined in
\cite[Appendix B]{PR}.
 
$\spadesuit$ In particular, we must study Fredholm operators 
in a ``skew-corner'' 
of our semifinite
von Neumann algebra $\mathcal N$. 
That is, if $P$ and $Q$ are projections in $\mathcal N$
(not necessarily infinite and not necessarily equivalent) we will extend
the notion of $\tau$-index and $\tau$-Fredholm to operators 
$T\in P\mathcal N Q$. If $\mathcal N$ is a factor, this is {\bf much easier}
and is done in Appendix A 
of \cite{Ph1}. We simply refer to them as 
$(P\cdot Q)$-Fredholm operators. Most results work in this 
setting; however the ploy used in \cite{Ph1} of invoking the 
existence of a partial isometry from $P$ to $Q$ to reduce to the case
$P\mathcal N P$ (solved in \cite{PR}) is not available. In fact, because of 
examples to which our version of the McKean-Singer Theorem applies, $P$ and
$Q$ are not generally equivalent. One notable result that is different in 
the nonfactor setting (even if $P=Q$) is that the set of $(P\cdot Q)$-Fredholm
operators with a given index {\bf is} open but is 
{\bf not} generally connected: information is lost when one fixes a 
trace to obtain a real-valued index. That the set of
$(P\cdot Q)$-Fredholm operators with a given index is open (and other facts) 
is very sensitive to the order in which the expected results are
proved. As the Fredholm alternative is {\bf not} available in
the $(P\cdot Q)$ setting, we take a novel approach and deduce many 
facts from the formula for the index of a product. We also study 
unbounded operators affiliated to a ``skew-corner''.
$\spadesuit$

\begin{not*}
If $T$ is an operator in the von Neumann algebra $\mathcal N$ (or $T$ is closed
and affiliated to $\mathcal N$) then we let
$R_T$ and $N_T$ be the projections on the closure of the range of $T$ and the 
kernel of $T$, respectively. If $T\in P\mathcal N Q,$ (or $T$ is closed and
affiliated to $P\mathcal N Q$) then we will denote
the projection on $ker_Q(T)=ker(T_{|_{Q(\HH)}})=ker(T)\cap Q(\HH)$ by $N^Q_T$ 
and observe that $N^Q_T=QN_T=N_TQ\leq Q$ while $R_T\leq P.$
\end{not*}

\begin{defn}
With the usual assumptionis on
$\mathcal N$ let $P$ and $Q$ be projections
({\bf not} necessarily infinite, or equivalent) in $\mathcal N$, and let 
$T\in P\mathcal N Q.$  
Then $T$ is called {\bf $(P\cdot Q)$-Fredholm} if and only if \vspace{.05in}\\
(1) $\tau(N_T^Q) < \infty$, and $\tau(N_{T^*}^P) < \infty$, and\vspace{.05in}\\
(2) There exists a $\tau$-finite projection $E\leq P$ with $range(P-E)\subseteq
range(T).$\vspace{.05in}\\
If $T$ is $(P\cdot Q)$-Fredholm then the {\bf $(P\cdot Q)$-index} of $T$ is
$$Ind(T)=\tau(N_T^Q)-\tau(N_{T^*}^P).$$
\end{defn}

\begin{lemma}\label{key}
With the usual assumptions on $\mathcal N$, let $T\in P\mathcal N Q.$ Then,
\vspace{.05in}\\
(1) With $P_1=R_T$ and $Q_1=Q-N^Q_T=supp(T)=R_{T^*}$, we have that $T$ is 
$(P\cdot Q)$-Fredholm if and only if $T$ is $(P_1\cdot Q_1)$-Fredholm
and in this case, the $(P_1\cdot Q_1)$-Index of $T$ is $0$, while the 
$(P\cdot Q)$-Index of $T$ is $\tau(Q-Q_1) - \tau(P-P_1)$.\vspace{.05in}\\
(2) If $T$ is $(P\cdot Q)$-Fredholm, then $T^*$ is $(Q\cdot P)$-Fredholm and
$Ind(T^*)=-Ind(T)$. If $T=V|T|$ is the polar decomposition,
then $V$ is $(P\cdot Q)$-Fredholm with $Ind(V)=Ind(T)$ and $|T|$ is 
$(Q\cdot Q)$-Fredholm of index $0$.\vspace{.05in}\\
(3) If $T=V|T|$ is $(P\cdot Q)$-Fredholm, then there exists a spectral 
projection $Q_0\leq Q$ for $|T|$ so that $\tau(Q-Q_0)<\infty$, and
$P_0=VQ_0V^*$ satisfies: $\tau(P-P_0)<\infty$, $P_0(\HH)=range(TQ_0)
\subset range(T)$, $Q_0(\HH)\subset range(T^*)$, 
$TQ_0=P_0TQ_0:Q_0(\HH)\to P_0(\HH)$ and $T^*P_0=Q_0T^*P_0:P_0(\HH)\to Q_0(\HH)$
are invertible as bounded linear operators.\vspace{.05in}\\ 
(4) The set of all $(P\cdot Q)$-Fredholm operators in $P\mathcal N Q$ is open 
in the norm topology.
\end{lemma}

\begin{proof}
(1) is straightforward, noting that $Q_1=1-N_T=R_{T^*}=supp(T)$.

(2) In the notation of part (1), $VV^*=P_1$ and $V^*V=Q_1$ so that $V$ is
$(P\cdot Q)$-Fredholm with $Ind(V)=Ind(T)$. Since both $T^*$ and $|T|$
have $\tau$-finite kernel and cokernel, it suffices to observe that if
$\tilde P\leq P$ is $\tau$-cofinite in $P$ and $\tilde P(\mathcal H)\subseteq
T(\mathcal H)$ then $\tilde Q:=V^*\tilde PV$ is $\tau$-cofinite in $Q$ and
satisfies $\tilde Q(\HH)\subseteq T^*(\HH)=|T|(\HH)$. The index statements 
are clear.

(3) By part (1), we can assume that $P=R_T$ and 
$Q=R_{T^*}=supp(T)$. Now $|T|\geq 0$ is 1:1 and
$\tau$-Fredholm in $Q\mathcal N Q$.
As $|T|$ is invertible modulo  
$\mathcal{K}_{Q\mathcal N Q}$ by Theorem B1 of \cite{PR}, the argument of 
Lemma 3.7
of \cite{CP0} shows that there exists a spectral 
projection $Q_0\leq Q$ for $|T|$ with $\tau(Q-Q_0)<\infty$ and $|T|Q_0$ is
bounded below on $Q_0(\HH)$. Let $P_0=VQ_0V^*$ this satisfies: 
$\tau(P-P_0)<\infty$. Now,
$TQ_0=V|T|Q_0=\cdots=P_0TQ_0,$ and
similarly, $T^*P_0=\cdots =Q_0T^*P_0.$ Since,
$TQ_0(\HH)=V|T|Q_0(\HH)=V(Q_0(\HH))=P_0(\HH),$ we see  
$P_0(\HH)=range(TQ_0)\subset range(T)$ and 
$TQ_0=P_0TQ_0:Q_0(\HH)\to P_0(\HH)$ is invertible as a bounded operator. 
The remaining bits are similar.

(4) Using (1), we have that $T$ is $(P_1\cdot Q_1)$-Fredholm of index $0$
and $V$ is a partial isometry in $\mathcal N$ with $VV^*=P_1$ and $V^*V=Q_1$.
By part (3) choose $Q_0\leq Q_1$ such that $\tau(Q_1-Q_0)<\infty$ and 
$P_0=VQ_0 V^*$ so that $\tau(P_1-P_0)<\infty$ satisfies 
$P_0(\HH)\subset range(T)$, and 
$TQ_0=P_0TQ_0$ is invertible as a bounded operator from $Q_0(\HH)$ 
to $P_0(\HH)$.
In particular, there exists $c>0$ so that for all $x\in\HH$:
$$\n TQ_0x\n=\n P_0TQ_0x\n\geq c\n Q_0x\n\;\;\&\;\;
\n T^*P_0x\n=\n Q_0T^*P_0x\n\geq c\n P_0x\n.$$
So if $A\in P\mathcal N Q$ and $\n T-A\n<c/3,$ then for all $x\in\HH$:
$$\n AQ_0x\n\geq 2c/3\n Q_0x\n\;\;\&\;\;
\n A^*P_0x\n\geq 2c/3\n P_0x\n.$$ 
Now clearly, $TQ_0$ and $T^*P_0$ have closed ranges $P_0$ and $Q_0,$
respectively. Let $\tilde{P}_0$ and $\tilde{Q}_0$ be the closed
ranges of $AQ_0$ and $A^*P_0$, respectively. Now, if $y\in P_0(\HH)$ is a unit 
vector, then $y=TQ_0x$ and $\n Q_0x\n\leq 1/c$. Letting
$y_1=AQ_0x\in \tilde{P}_0(\HH)$, we have $\n y-y_1\n\leq (c/3)(1/c)=1/3.$
Similarly, if $z\in \tilde{P}_0(\HH)$ is a unit vector, we find 
$z_1\in P_0(\HH)$ with $\n z_1-z\n\leq (c/3)(3/2c)=1/2.$
On concludes that $\n P_0-\tilde{P}_0\n\leq 1/3+1/2 <1,$ and so
$P_0$ and $\tilde{P}_0$ are unitarily equivalent by a unitary in
$P\mathcal N P$ that fixes $P$. Hence, $P-\tilde{P}_0$ is $\tau$-cofinite
and not only is $\tilde{P}_0(\HH)\subset range(A)$, but also
$N_{A^*}^P=P-R_A\leq (P-\tilde{P}_0)$ is $\tau$-finite. Similarly, 
$N_A^Q\leq (Q-\tilde{Q}_0)$
is $\tau$-finite and $A$ is $(P\cdot Q)$-Fredholm. 
\end{proof}

\begin{defn}
If $T\in P\mathcal N Q$, then a {\bf parametrix} for $T$ is an operator
$S\in Q\mathcal N P$ satisfying $ST=Q+k_1$ and $TS=P+k_2$ where 
$k_1\in\mathcal{K}_{Q\mathcal N Q}$ and $k_2\in\mathcal{K}_{P\mathcal N P}.$
\end{defn}

\begin{lemma}\label{para}
With the usual assumptions on $\mathcal N$, then $T\in P\mathcal N Q$ is
$(P\cdot Q)$-Fredholm if and only if $T$ has a parametrix $S\in Q\mathcal N P$.
Moreover, any such parametrix is $(Q\cdot P)$-Fredholm.
\end{lemma}
\begin{proof}
Let $S$ be a parametrix for $T$. Then $TS=P+k_2$ is Fredholm in $P\mathcal N P$
by Appendix B of \cite{PR}. Hence there exists a projection $P_1\leq P$
with $\tau(P-P_1)<\infty$ and $P_1(\HH)\subset range(TS)\subset range(T).$
So, $N_{T^*}^P = P-R_T\leq P-P_1$ is $\tau$-finite. On the other hand,
$T^*S^*=(ST)^*=Q+k_1^*$ is Fredholm in $Q\mathcal N Q$ again by Appendix B
of \cite{PR} and so by the same argument $N_T^Q$ is also $\tau$-finite.
That is, $T$ is $(P\cdot Q)$-Fredholm and similarly $S$ is 
$(Q\cdot P)$-Fredholm.

Now suppose that $T$ is $(P\cdot Q)$-Fredholm. By part (3) of Lemma \ref{key}, 
there
exist projections $Q_0$ and $P_0$ which are $\tau$-cofinite in $Q$ and $P$
respectively so that $TQ_0=P_0TQ_0:Q_0(\HH)\to P_0(\HH)$ is invertible as
a bounded linear operator. Let $S$ be its inverse. Then $S\in\mathcal N$
so that $S=Q_0SP_0\in Q\mathcal N P$, and $STQ_0=Q_0$ and $TQ_0S=P_0$.
Finally, 
$$ST=STQ_0+ST(Q-Q_0)=Q_0+k=Q+k_1\;\;and\;\;TS=TQ_0S=P_0=P+k_2,$$
where $k_1\in\mathcal{K}_{Q\mathcal N Q}$ and 
$k_2\in\mathcal{K}_{P\mathcal N P}.$ That is, $S$ is a parametrix for $T$.
\end{proof}

\begin{lemma}\label{compact}
We retain the usual assumptions on $\mathcal N$.\vspace{.05in}\\ 
(1) Let $T\in P\mathcal N Q$ be $(P\cdot Q)$-Fredholm. If 
$k\in P\mathcal{K}_{\mathcal N}Q$
then $T+k$ is also $(P\cdot Q)$-Fredholm.\vspace{.05in}\\
(2) If $T\in P\mathcal N Q$ is $(P\cdot Q)$-Fredholm and $S\in G\mathcal N P$
is $(G\cdot P)$-Fredholm, then $ST$ is $(G\cdot Q)$-Fredholm.
\end{lemma}

\begin{proof}
One checks that if $S$ is a parametrix for $T$ then $S$ is also 
a parametrix for $T+k$ and that if $T_1$ is a parametrix for $T$ and
$S_1$ is a parametrix for $S$, then $T_1S_1$ is a parametrix for $ST$.
\end{proof}
 
\begin{prop}\label{productform}
Let $G,P,Q$ be projections in $\mathcal N$ (with 
trace $\tau$) and let $T\in P\mathcal N Q$ be $(P\cdot Q)$-Fredholm and 
$S\in G\mathcal N P$ be $(G\cdot P)$-Fredholm, respectively. Then,
$ST$ is $(G\cdot Q)$-Fredholm and
$$Ind(ST)=Ind(S)+Ind(T).$$
\end{prop}
We follow Breuer in \cite{B2} indicating the changes needed in this generality. 
Before proving the proposition we require a Lemma. 
\begin{lemma}(Cf. Lemma 1 of \cite{B2}) 
With the hypotheses of the Proposition:
$$N_{ST}^Q - N_T^Q\sim inf(R_T,N_S^P).$$
\end{lemma}

\begin{proof}
We follow Breuer's arguments replacing $ker(T)$ with 
$ker_Q(T)=ker(T)\cap Q(\HH)$; $ker(ST)$ with $ker_Q(ST)$; and $ker(S)$ 
with $ker_P(S)$. Noting $N_T^Q=QN_T=N_T Q$ and similar identities, we read 
Breuer until we choose projections 
$E_1\leq E_2\leq E_3\leq...$ as in Lemma 13 of
\cite{B1} satisfying each $(P-E_n)$ is $\tau$-finite, 
$E_n(\HH)\subset range(T)$, and $sup\{E_n|\;  n=1,2,...\}=R_T$.
We continue reading carefully, replacing $1$ with $P$ at crucial points.
Finally, we get the conclusion from:\vspace{.05in}\\
\hspace*{1.4in}$N_{ST}^Q - N_T^Q=R_{(N_{ST}^Q-N_T^Q)T^*}\sim 
R_{T(N_{ST}^Q-N_T^Q)}
=inf(R_T,N_S^P).$\end{proof}

\begin{proof} (Of the Proposition):\\
Now, $S^*$, $T^*$, $ST$, and $(ST)^*=T^*S^*$ are all Fredholm, 
and the above lemma implies:
$$N_{ST}^Q-N_T^Q \sim inf(R_T,N_S^P)\;\;
and
\;\;N_{(ST)^*}^G-N_{S^*}^G \sim inf(R_{S^*},N_{T^*}^P).$$
The projections on the RHS of the two similarities are in $P\mathcal N P$, 
and so by \cite[Cor. 1, p.216]{Dix}:
$$N_S^P-inf(P-N_{T^*}^P,N_S^P)\sim N_{T^*}^P- inf(P-N_S^P,N_{T^*}^P).$$
Since $P-N_{T^*}^P=R_T$ and $P-N_S^P=R_{S^*},$ we get:
$$N_S^P-inf(R_T,N_S^P)\sim N_{T^*}^P- inf(R_{S^*},N_{T^*}^P).$$
Using these similarities we calculate:
$$Ind(ST)=\tau(N_{ST}^Q)-\tau(N_{(ST)^*}^G)=
\tau(N_{ST}^Q-N_T^Q)-\tau(N_{(ST)^*}^G-N_{S^*}^G)+
\tau(N_T^Q)-\tau(N_{S^*}^G)$$
\hspace*{1in}$=\cdots=\tau(N_S^P)-\tau(N_{S^*}^G) + \tau(N_T^Q)-\tau(N_{T^*}^P)=
Ind(S) + Ind(T).$
\end{proof}

\begin{cor}\label{invariance} (Invariance properties of the $(P\cdot Q)$-Index)
Let $T\in P\mathcal N Q.$
\vspace{.05in}\\
(1) If $T$ is $(P\cdot Q)$-Fredholm then there exists $\delta>0$ so that
if $S\in P\mathcal N Q$ and $\n T-S\n<\delta$ then $S$ is 
$(P\cdot Q)$-Fredholm and $Ind(S)=Ind(T).$\\
(2) If $T$ is $(P\cdot Q)$-Fredholm and
$k\in P\mathcal K_{\mathcal N}Q$ then $T+k$ is
$(P\cdot Q)$-Fredholm and\\ $Ind(T+k)=Ind(T).$
\end{cor}

\begin{proof}
(1)  By the Proposition and part (2) of Lemma \ref{key},
$TT^*$ is Fredholm of index $0$ in $P\mathcal N P$. So by 
Corollary B2 of \cite{PR} there exists $\epsilon_1>0$ so that if 
$A\in P\mathcal N P$
satisfies $\n A-TT^*\n<\epsilon_1$ then $A$ is Fredholm of
index $0$. Moreover,
by part (4) of Lemma \ref{key} there exists $\epsilon_2>0$ so that the ball of
radius $\epsilon_2$ about $T$ in $P\mathcal N Q$ is contained in the
$(P\cdot Q)$-Fredholms. Let 
$\delta = min \{\epsilon_2,\epsilon_1/||T||\}$. Then if $S\in P\mathcal N Q$ 
and $\n T-S\n<\delta$ then $S$ is $(P\cdot Q)$-Fredholm and 
$\n ST^*-TT^*\n<\epsilon_1$ so that $ST^*$ is $(P\cdot P)$-Fredholm of index 
$0$. By the Proposition and part (2) of Lemma \ref{key}:
$$0=Ind(ST^*)=Ind(S)-Ind(T).$$
(2) This is similar to part (1) but uses Lemma \ref{compact} part (1) in place 
of Lemma \ref{key} part (4).
\end{proof}

$\spadesuit$ In \cite{Ph1} {\bf spectral flow} is defined in a
semifinite {\bf factor} using the index of
Breuer-Fredholm operators in a skew-corner $P\mathcal NQ$ (in particular the
operator $PQ$) and uses the product theorem for the index and other
standard properties. The non-factor case for Toeplitz operators ($P=Q$) 
is covered in \cite{PR} but the more subtle ``skew-corner'' case
has not appeared in the literature.
This section enables one to
extend \cite{Ph1} to the nonfactor setting where it was needed for
\cite{CP2}, \cite{CPS2} and \cite{CPRS2}.
For use in the present paper we generalise some of these results to 
closed, densely defined operators affiliated to $P\mathcal N Q$ by 
studying the map $T\mapsto T(1+|T|^2)^{-1/2}$.$\spadesuit$

\begin{defn} A closed, densely defined operator $T$ affiliated to 
$P\mathcal N Q$
is {\bf $(P\cdot Q)$-Fredholm} if\vspace{.05in}\\
(1) $\tau(N_T^Q) < \infty$, and $\tau(N_{T^*}^P) < \infty$, and\vspace{.05in}\\
(2) There exists a $\tau$-finite projection $E\leq P$ with $range(P-E)\subset
range(T).$\vspace{.05in}\\
If $T$ is $(P\cdot Q)$-Fredholm then the {\bf $(P\cdot Q)$-index} of $T$ is:
$Ind(T)=\tau(N_T^Q)-\tau(N_{T^*}^P).$
\end{defn}

{\bf Remark} Using the equalities:
$range(1+|T|^2)^{-1/2}=dom(1+|T|^2)^{1/2}=dom(|T|)=dom(T)$
one can show that:\; 
$range(T)=range(T(1+|T|^2)^{-1/2})$;\; 
$ker(T)=ker(T(1+|T|^2)^{-1/2})$\;and 
\;$ker(T^*)=ker([T(1+|T|^2)^{-1/2}]^*).$
A little more thought completes the following:

\begin{prop}\label{transform}{\bf (Index)}
If $T$ is a closed, densely defined operator affiliated to $P\mathcal N Q$, 
then
$T$ is $(P\cdot Q)$-Fredholm if and only if the operator 
$T(1+|T|^2)^{-1/2}$ is $(P\cdot Q)$-Fredholm in $P\mathcal N Q$. In this case, 
$$Ind(T)=Ind(T(1+|T|^2)^{-1/2}).$$
\end{prop}

\begin{prop}{\bf (Continuity)} If $T$ is a closed, densely defined operator 
affiliated to 
$P\mathcal N Q$, and $A\in P\mathcal N Q$ then $T+A$ is also 
closed, densely defined, and affiliated to $P\mathcal N Q$ and
$$\n T(1+|T|^2)^{-1/2} - (T+A)(1+|T+A|^2)^{-1/2}\n\leq\n A\n.$$
\end{prop}

\begin{proof}
We define the following self-adjoint operators:
$$D=\left(\begin{array}{cc} 0 &  T^*  \\
                             T &  0
  \end{array} \right)\;\;\;and\;\;\;
  B=\left(\begin{array}{cc} 0 &  A^*  \\
                             A &  0
  \end{array} \right).$$
Then, $D$ is affiliated to $M_2(\mathcal N)$ and $B\in M_2(\mathcal N)$.
By \cite[Theorem 8, Appendix A]{CP1}, we have:
$$\n D(1+D^2)^{-1/2} - (D+B)(1+(D+B)^2)^{-1/2}\n\leq\n B\n.$$
A little calculation yields:
$$\n T(1+|T|^2)^{-1/2} - (T+A)(1+|T+A|^2)^{-1/2}\n$$
\hspace*{1.6in}$\leq\n D(1+D^2)^{-1/2} - (D+B)(1+(D+B)^2)^{-1/2}\n
\leq\n B\n =\n A\n.$
\end{proof}

\begin{cor}{\bf (Index continuity)}
If $T$ is  affiliated to $P\mathcal N Q$ and
$T$ is $(P\cdot Q)$-Fredholm then there exists $\epsilon>0$ so that if 
$A\in P\mathcal N Q$ and $\n A\n<\epsilon$, then $T+A$ is 
$(P\cdot Q)$-Fredholm and
$$Ind(T+A)=Ind(T).$$
\end{cor}

\begin{prop}{\bf (Compact perturbation)}
Let $T$ be {\bf any} closed, densely defined operator affiliated to 
$P\mathcal N Q.$
\vspace{.05in}\\
(1)\;\;If $k\in P\mathcal K_{\mathcal N}Q$, then the difference 
$T(1+|T|^2)^{-1/2} - (T+k)(1+|T+k|^2)^{-1/2}$ is in
$P\mathcal K_{\mathcal N}Q$!\vspace{.05in}\\
(2)\;\;If $T$ is $(P\cdot Q)$-Fredholm then for all 
$k\in P\mathcal K_{\mathcal N}Q$, $T+k$ is $(P\cdot Q)$-Fredholm
and $$Ind(T+k)=Ind(T).$$

\end{prop}

\begin{proof}
We prove the surprisingly subtle (and rather surprising!)
first statement, since part (2) is an immediate 
corollary by Proposition \ref{transform} and Corollary \ref{invariance}.
By the $2\times 2$ matrix trick, we can assume that
$T$ and $k$ are self-adjoint and that $P=Q=1.$ By the resolvent equation:
$$(T+i1)^{-1}-(T+k+i1)^{-1}=(T+i1)^{-1}k(T+k+i1)^{-1}\in 
\mathcal K_{\mathcal N}.$$
However the identity, $(T+i1)^{-1}=T(1+T^2)^{-1}-i(1+T^2)^{-1}$ and the 
corresponding identity for $(T+k+i1)^{-1}$ imply that:
$$T(1+T^2)^{-1}-(T+k)(1+(T+k)^2)^{-1}\in \mathcal K_{\mathcal N}$$ since
this difference is the self-adjoint part of an element in the $C^*$-algebra
$\mathcal K_{\mathcal N}.$ Now for $\mu>0$ real we can replace $T$ with
$\mu T$ and $(T+k)$ with $\mu(T+k)$ and get:
$$\mu\left\{\mu T(1+(\mu T)^2)^{-1}-\mu(T+k)(1+(\mu(T+k))^2)^{-1}\right\}
\in\mathcal K_{\mathcal N}.$$
For any real $\lambda\geq 0$, let $\mu = (1+\lambda)^{-1/2}$ and a little 
calculation yields:
$$T(1+T^2+\lambda)^{-1}-(T+k)(1+(T+k)^2+\lambda)^{-1}
\in\mathcal K_{\mathcal N}.$$
By \cite[Lemma 6, Appendix A]{CP1} we have the estimate:
$$\n T(1+T^2+\lambda)^{-1}-(T+k)(1+(T+k)^2+\lambda)^{-1}\n\leq\frac{\n k\n}
{1+\lambda},$$
and hence the following integral converges absolutely in operator norm
to an element in the $C^*$-algebra $\mathcal K_{\mathcal N}:$
$$\frac{1}{\pi}\int_0^\infty \lambda^{-1/2}\left(
T(1+T^2+\lambda)^{-1}-(T+k)(1+(T+k)^2+\lambda)^{-1}\right)d\lambda.$$
If we call this element $k_0$, then by \cite[Lemma 4, Appendix A]{CP1},
we have for all $\xi\in dom(T)=dom(T+k)$ that the following integrals 
{\bf converge in} ${\mathbf \HH}$ and:
\bean
&&T(1+T^2)^{-1/2}(\xi)-(T+k)(1+(T+k)^2)^{-1/2}(\xi)\nno
&=&\frac{1}{\pi}\int_0^\infty \lambda^{-1/2}
T(1+T^2+\lambda)^{-1}(\xi)d\lambda-\frac{1}{\pi}\int_0^\infty \lambda^{-1/2}
(T+k)(1+(T+k)^2+\lambda)^{-1}(\xi)d\lambda\nno
&=&\frac{1}{\pi}\int_0^\infty \lambda^{-1/2}\left(
T(1+T^2+\lambda)^{-1}-(T+k)(1+(T+k)^2+\lambda)^{-1}\right)(\xi) d\lambda=
k_0(\xi).
\eean
As both side of this equation are bounded operators, we have:\vspace{.05in}\\
\hspace*{1.5in}
$T(1+T^2)^{-1/2}-(T+k)(1+(T+k)^2)^{-1/2}=k_0\in\mathcal K_{\mathcal N}.$
\end{proof}

\begin{defn}
For many geometric examples, the following is a useful notion. If $T$ is a 
closed, densely defined, unbounded operator
affiliated to $P\mathcal N Q$ then a {\bf parametrix} for $T$ is a bounded
everywhere defined operator $S\in Q\mathcal N P$ so that:\vspace{.05in}\\
(1) $\overline{TS}=P+k_1$ for $k_1\in P\mathcal K_{\mathcal N}P$,
\vspace{.05in}\\
(2) $\overline{ST}=Q+k_2$ for $k_2\in Q\mathcal K_{\mathcal N}Q.$
\vspace{.05in}\\
{\bf Note}, since $T$ is closed and $S$ is bounded, $TS=\overline{TS}$ is
everywhere defined and bounded by $(1)$. For example, if $D$ is an
unbounded self-adjoint operator and $(1+D^2)^{-1}\in {\mathcal K}_{\mathcal N}$
then $D(1+D^2)^{-1}$ is a parametrix for $D$ since 
$\overline{D(1+D^2)^{-1}D}=D^2(1+D^2)^{-1}=1-(1+D^2)^{-1}$.
\end{defn}

\begin{lemma}
If $T$ is a closed, densely defined, unbounded operator
affiliated to $P\mathcal N Q$ then $T$ has a parametrix if and only if $T$ is
$(P\cdot Q)$-Fredholm.
\end{lemma}
\begin{proof}
If $S$ is a parametrix for $T$ then by (1) $TS$ is everywhere defined and 
Fredholm in $P\mathcal N P$. So there exists a projection
$E\leq P$ with $\tau(E)<\infty$ and:
$range(P-E)\subset range(TS)\subset range(T).$
In particular, this implies (since $TS$ is bounded) that
$N_{(TS)^*}^P$ is $\tau$-finite. But $S^*T^*\subseteq (TS)^*$ and so
$N^P_{T^*}\leq N^P_{(TS)^*}.$ 
That is, $\tau(N_{T^*}^P)<\infty.$
Now, $\overline{ST}=Q+k_2$ is $(Q\cdot Q)$-Fredholm and so has a $\tau$-finite
$Q$-kernel. But $N^Q_T\leq N^Q_{\overline{ST}}.$
That is, $\tau(N_T^Q)<\infty$ and $T$ is $(P\cdot Q)$-Fredholm.\\
\hspace*{.2in}If $T=V|T|$ is $(P\cdot Q)$-Fredholm then 
$|T|(1+|T|^2)^{-1/2}$ is bounded and
$(Q\cdot Q)$-Fredholm and so has a parametrix $S$ which we can take to be 
a function of $|T|(1+|T|^2)^{-1/2}$. Thus $S$ commutes
with $(1+|T|^2)^{-1/2}$. One then checks that $(1+|T|^2)^{-1/2}SV^*$ is a 
parametrix for $T$.
\end{proof}

{\bf Remark} In general a parametrix for a genuinely unbounded Fredholm operator
is {\bf not} Fredholm as its range cannot contain the range of a cofinite
projection.

\begin{thm}
(McKean-Singer) Let $D$ be an unbounded self-adjoint operator affiliated
to the semifinite von Neumann algebra $\mathcal N$ (with faithful normal 
semifinite trace $\tau$). Let $\gamma$ be a self-adjoint unitary in $\mathcal N$
which anticommutes with $D$. Finally, let $f$ be a continuous even function 
on $\bf{R}$ with $f(0)\neq 0$ and
$f(D)$ trace-class. Let $D^+=P^{\perp}D P$ where $P=(\gamma +1)/2$
and $P^{\perp}=1-P$. Then
as an operator affiliated to $P^{\perp}\mathcal N P$, $D^+$ is 
$(P^{\perp}\cdot P)$-Fredholm and\vspace{.05in}\\ 
\hspace*{2.5in}$Ind(D^+)=\frac{1}{f(0)}\tau\left(\gamma f(D)\right).$
\end{thm}

\begin{proof}
Let $D^-=PD P^{\perp}$. Since $\{D,\gamma\}=0$, we see that relative to the
decomposition $1=P\oplus P^{\perp}$:
$$\gamma =\left(\begin{array}{cc} 1 & 0\\
                                  0 & -1
  \end{array} \right),\;\;\;\\
  D=\left(\begin{array}{cc} 0   &  D^-\\
                          D^+ &  0
  \end{array} \right),\;\;\;\\
  D^2=\left(\begin{array}{cc} D^-D^+ & 0\\
                                0    & D^+D^-
  \end{array} \right),\;\;\;\\
  |D|=\left(\begin{array}{cc} |D^+| & 0\\
                                0    & |D^-|
  \end{array} \right).$$ 
  
We have already observed that $D(1+D^2)^{-1}$ is a parametrix for $D$. But, 
then:
$$D(1+D^2)^{-1}=\left(\begin{array}{cc}  0 & D^-(P^{\perp}+|D^-|^2)^{-1}\\ 
                       D^+(P+|D^+|^2)^{-1} & 0
  \end{array}\right).$$ 
Hence $D^-(P^{\perp}+|D^-|^2)^{-1}$ is a parametrix for $D^+$ and so $D^+$ is 
$(P^{\perp}\cdot P)$-Fredholm.
   Let $D^+=V|D^+|$ be the polar decomposition of $D^+$ so that
$D^-=D^{+*}=|D^+|V^*$. Then $V\in\mathcal N$ is a partial isometry
with initial space $P_1=V^*V=supp(D^+)\leq P$ and final space
$Q_1=VV^*=range(D^+)^{-}=supp(D^-)\leq P^{\perp}.$ Then, 
$ker(D^+)=P_0(\mathcal H)$ 
as an operator on $P(\mathcal H)$ where $P_0=P-P_1.$   
Similarly, $coker(D^+)=ker(D^-)=Q_0(\mathcal H)$ where
$Q_0=P^{\perp}-Q_1$.  

Now, $|D^+|^2=D^-D^+=D^-D^{-*}=V^*|D^-|^2V$ so that $|D^+|=V^*|D^-|V$ 
and if $g$
is any bounded continuous function then,  $g(|D^+|_{|_{P_1(\mathcal H)}})=
V^*g(|D^-|_{|_{Q_1(\mathcal H)}})V.$  But, as an operators on 
$P(\mathcal H)$, and respectively, $P^{\perp}(\mathcal H)$ we have:
\bean g(|D^+|)&=& P_1g(|D^+|)P_1\oplus g(0)P_0=
g(|D^+|_{|_{P_1(\mathcal H)}})\oplus g(0)P_0\;\;and\\
g(|D^-|)&=& Q_1g(|D^-|)Q_1\oplus g(0)Q_0=
g(|D^-|_{|_{Q_1(\mathcal H)}})\oplus g(0)Q_0.\eean
Finally, since $f$ is even, we have $f(D)=f(|D|)$ and so:
$$\gamma f(D)=\left(\begin{array}{cc} f(|D^+|) & 0\\
                                0    & -f(|D^-|)
  \end{array} \right) = \left(\begin{array}{cc} f(|D^+|_{|_{P_1(\mathcal H)}})
  \oplus f(0)P_0) & 0\\
                                0    & -f(|D^-|_{|_{Q_1(\mathcal H)}})\oplus
         -f(0)Q_0                       
  \end{array} \right)$$
$$= \left(\begin{array}{cc} V^*f(|D^-|_{|_{Q_1(\mathcal H)}})V
  \oplus f(0)P_0) & 0\\
                                0    & -f(|D^-|_{|_{Q_1(\mathcal H)}})\oplus
         -f(0)Q_0 \end{array} \right).$$
Hence,\hspace{1in} $\tau(\gamma f(D)) = f(0)\tau(P_0)-f(0)\tau(Q_0) = 
f(0)Ind(D^+).$                  
\end{proof}

\begin{cor} Let $(\A,\HH,\D)$ be an even spectral triple with grading
$\gamma$, $(1+\D^2)^{-1/2}\in \mathcal{L}^n(\mathcal N)$ and $p\in\A$,
a projection. Then, relative to the decomposition afforded by 
$\gamma$ as above, we have: 
$$p=\left(\begin{array}{cc} p^{+} & 0\\
                                0    &  p^{-}
  \end{array} \right),\;\;where\;\; p^{+}=PpP=Pp\;\;and\;\; 
p^{-}=P^{\perp}pP^{\perp}=pP^{\perp}.$$
So, $p\D^+ p=pP^{\perp}\D Pp=p^{-}\D p^+$ is an operator affiliated to 
$p^{-}\mathcal N p^{+}$ we have that 
$p^{-}\D^{+}p^{+}$ is $(p^{-}\cdot p^{+})$-Fredholm and for
any fixed $a\geq 0$ its  $\mathbf {(p^{-}\cdot p^{+})}${\bf -index} is given by:
$$Ind(p\D^+ p)=Ind(p^{-}\D^{+}p^{+})=(1+a)^{n/2}
\tau\left(\gamma p\left(p+a+(p\D p)^2\right)^{-n/2}\right).$$ 
\end{cor}

\begin{proof}
In the above version of the McKean-Singer theorem, we replace $\A$ 
with $p\A p$ which is a unital
subalgebra of the semifinite von Neumann algebra $p\mathcal N p$. Moreover,
the operator $p\D p$ is self-adjoint and affiliated to $p\mathcal N p$,
and $p\gamma$ is a grading in $p\mathcal N p$. One easily checks that
$$(p\D p)^+ = p^{-}\D^+p^{+}.$$
Letting $f(x)=(1+a+x^2)^{-n/2}$, we can apply the McKean-Singer theorem
once we show that $(p+a+(p\D p)^2)^{-1/2}\in \mathcal{L}^n(p\mathcal N p).$
It suffices to do this for $a=0$ since $$(p+a+(p\D p)^2)^{-1/2}\leq 
(p+(p\D p)^2)^{-1/2}.$$
This is a careful calculation:
\bean
&&p(1+\D^2)^{-1}p-(p+(p\D p)^2)^{-1}\nno
&=& p[(1+\D^2)^{-1} - \{(p+(p\D p)^2) + (1-p)\}^{-1}]p\nno
&=& p(1+\D^2)^{-1}\left((p+(p\D p)^2)+(1-p)-(1+\D^2)\right)
\{(p+(p\D p)^2)+(1-p)\}^{-1}p\nno
&=& p(1+\D^2)^{-1}\left((p\D p)^2-\D^2\right)p(p+(p\D p)^2)^{-1}p\nno
&=& p(1+\D^2)^{-1}\left([p,\D]p\D p+\D p[\D,p]+\D[p,\D]\right)
p(p+(p\D p)^2)^{-1}p\nno
&=& p(1+\D^2)^{-1}[p,\D]p\D p(p+(p\D p)^2)^{-1}p\nno
&& \hspace{.6in} + p(1+\D^2)^{-1}\D p[\D,p]p(p+(p\D p)^2)^{-1}p\nno
&& \hspace{1.2in} + p(1+\D^2)^{-1}\D[p,\D]p(p+(p\D p)^2)^{-1}p
\eean
Now, since $|\D(1+\D^2)^{-1}|\leq (1+\D^2)^{-1/2}$, we have that three terms in 
the last lines are in $\mathcal L^{n/2}$, $\mathcal L^{n}$, and 
$\mathcal L^{n}$, respectively, and so their sum is in $\mathcal L^{n}.$
Since, $p(1+\D^2)^{-1}p \in \mathcal L^{n/2}$, we see from the first line
in the displayed equations that 
$(p+(p\D p)^2)^{-1}$ is in $\mathcal L^{n}.$\\
\hspace*{.2in}Now, armed with this new information, we look at the three 
terms in the last 
line again, and see that they are in $\mathcal L^{n/2}$,
$\mathcal L^{n}\cdot\mathcal L^{n}$, and $\mathcal L^{n}\cdot\mathcal L^{n}$,
respectively, and so their sum is in $\mathcal L^{n/2}$. Thus,
$(p+(p\D p)^2)^{-1}$ is, in fact, in $\mathcal L^{n/2}$: in other words,
$(p+(p\D p)^2)^{-1/2}$ is in $\mathcal L^{n}$ as claimed. 
\end{proof}

$\spadesuit$ From now on, we follow convention and denote the above index
by $Ind(p\D^+p)$; effectively disguising the fact that $p\D^+ p$
is, in fact, Fredholm relative to the ``skew-corner,'' $p^{-}\mathcal N p^{+}$.
$\spadesuit$

{\bf Remark:} The ideal $\mathcal L^{n}(\mathcal N)$ 
can be replaced by any symmetric ideal $\mathcal I\subset\mathcal K_{N}$ 
provided we use an even function $f$ satisfying 
$f(|T|)\in \mathcal L^1$ for all $T\in\mathcal I$. The formula then
becomes:\\
\hspace*{1.6in}$Ind(p\D^+ p)=(1/f(0))
\tau\left(\gamma pf\left((\left(p+(p\D p)^2\right)^{-1/2}\right)\right).$\\
In particular, if $(\A,\HH,\D)$ is $\theta$-summable, and $f(x)=e^{-tx^2}$, 
$t>0$, the formula becomes:\\
\hspace*{2.4in}$Ind(p\D^+ p)=\tau\left(\gamma pe^{-t(p\D p)^2}\right).$

\section{Statement of the Main Result}\label{mainresult}

We use the notation of \cite{CPRS2}.
Denote multi-indices by $(k_1,...,k_m)$, $k_i=0,1,2,...$, whose length 
$m$ will always be 
clear from the context and let $|k|=k_1+\cdots+k_m$. Define
\ben \alpha(k)=({k_1!k_2!\cdots k_m!(k_1+1)(k_1+k_2+2)\cdots(|k|+m)})^{-1}
\een
and $\s_{n,j}$ (the elementary symmetric functions of $\{1,...,n\}$) by
$\prod_{j=0}^{n-1}(z+j)=\sum_{j=1}^{n}z^j\s_{n,j}$.
If  $(\A,\HH,\D)$ is a $QC^\infty$ spectral triple and $T\in\cn$ 
then $T^{(n)}$ is the $n^{th}$ iterated commutator with $\D^2$, that is,  
$[\D^2, [\D^2,[\cdots,[\D^2,T]\cdots]]]$. 

We let $q=\inf\{k\in {\mathbf R}:\tau((1+\D^2)^{-k/2})<\infty\}$
be the {\bf spectral dimension} of $(\A,\HH,\D)$ and we assume it is 
{\bf isolated},
{\it ie}, 
for  $$b=a_0[\D,a_1]^{(k_1)}\cdots[\D,a_m]^{(k_m)}(1+\D^2)^{-m/2-|k|}$$
the zeta functions
\ben \zeta_b(z-(1-q)/2)=\tau(b(1+\D^2)^{-z+(1-q)/2})\een
have analytic continuations to a deleted neighbourhood of $z=(1-q)/2$.
As in \cite{CPRS2} we let
$\tau_j(b)=res_{z=(1-q)/2}(z-(1-q)/2)^j\zeta_b(z-(1-q)/2)$.
Our main result is:

\begin{thm}[Semifinite Even Local Index Theorem]\label{SFLIT} 
Let $(\A,\HH,\D)$ be an 
even $QC^\infty$ spectral triple with 
spectral dimension
$q\geq 1$.  Let $N=[\frac{q+1}{2}]$, where $[\cdot]$ denotes the integer part, 
and let $p\in\A$ be a self-adjoint projection. Then

1) \qquad\qquad\qquad $Ind(p\D^+p)= res_{r=(1-q)/2}\left( 
\sum_{m=0,even}^{2N} \phi_m^r(Ch_m(p))\right)$\\
where for $a_0,...,a_m\in\A$, $l=\{a+iv:v\in{\R}\}$, $0<a<1/2$, 
$R_s(\lambda)=(\lambda-(1+s^2+\D^2))^{-1}$ and $r>1/2$ we define 
$\phi_m^r(a_0,a_1,...,a_m)$ to be
\ben \frac{(m/2)!}{m!}\int_0^\infty 2^{m+1}s^m\tau\left(\gamma
\frac{1}{2\pi i}\int_l\lambda^{-q/2-r}a_0R_s(\lambda)[\D,a_1]R_s(\lambda)
\cdots[\D,a_m]R_s(\lambda)d\lambda\right)ds.\een
In particular the sum on the right hand side of $1)$ analytically continues to 
a deleted neighbourhood of $r=(1-q)/2$ with {\em at worst} a simple pole at 
$r=(1-q)/2$.
Moreover, the complex function-valued cochain $(\phi_m^r)_{m=0,even}^{2N}$ is 
a $(b,B)$ 
cocycle for $\A$ modulo functions holomorphic in a half-plane containing 
$r=(1-q)/2$.

2) The index, $Ind(p\D^+p)$ is also the residue of a sum of zeta functions:
\bean && res_{r=(1-q)/2} \Biggl(\sum_{m=0,even}^{2N}\sum_{|k|=0}^{2N-m}
\sum_{j=1}^{|k|+m/2}(-1)^{|k|+m/2}\alpha(k)
\frac{(m/2)!}{2m!}\s_{|k|+m/2,j}\times \Biggr.\nno
&& \quad\times\Biggl.(r-(1-q)/2)^j\tau\left(\gamma(2p-1)
[\D,p]^{(k_1)}[\D,p]^{(k_2)}\cdots[\D,p]^{(k_m)}(1+\D^2)^{-m/2-|k|-r+(1-q)/2}
\right)\Biggr),\eean
(for $m=0$ we replace $(2p-1)$ by $2p$).
In particular the sum of zeta functions on the right hand side analytically 
continues to a deleted neighbourhood of $r=(1-q)/2$ and has {\rm at worst} a 
simple pole at 
$r=(1-q)/2$.

3) If $(\A,\HH,\D)$ also has isolated spectral dimension then
\ben Ind(p\D^+p)=\sum_{m=0,even}^{2N} \phi_m(Ch_m(p))\een
where for $a_0,...,a_m\in\A$ we have $\phi_0(a_0)=
res_{r=(1-q)/2}\phi^r_0(a_0)=\tau_{-1}(\gamma a_0)$ and for $m\geq 2$
\bean \phi_m(a_0,...,a_m)&=&res_{r=(1-q)/2}\phi^r_m(a_0,...,a_m)=
\sum_{|k|=0}^{2N-m}(-1)^{|k|}
\alpha(k)\times\nno
&\times&\sum_{j=1}^{|k|+m/2}\s_{(|k|+m/2),j}\tau_{j-1}
\left(\gamma a_0[\D,a_1]^{(k_
1)}\cdots[\D,a_m]^{(k_m)}(1+\D^2)^{-|k|-m/2}\right),\eean
and $(\phi_m)_{m=0,even}^{2N}$ is a $(b,B)$ cocycle for $\A$. When $[q]=2n+1$ 
is odd, the term with $m=2N$ is zero, and 
for $m=0,2,...,2N-2$, all the top terms with $|k|=2N-m$ are zero. 
\end{thm}

\begin{cor}\label{lowdim}
For $1\leq q<2$, the statements in $3)$ of Theorem \ref{SFLIT} are true 
without the 
assumption of isolated dimension spectrum.
\end{cor}

\section{The Local Index Theorem in the Even Case}\label{evencase}

The main technical device
that improves the proof of the local index theorem of \cite{CM} for odd 
spectral triples  stems from  our 
use in \cite{CPRS2} of the resolvent cocycle
to reduce the hypotheses needed for the theorem and most importantly
to provide a simple proof that  
the (renormalised) residue cocycle of Connes-Moscovici is an index cocycle. 
We will see that these improvements also apply in the even case.

In this Section we will derive the formulae for the index
appearing in parts $2)$ and $3)$ of Theorem \ref{SFLIT}.
The exposition is broken down into six subsections. 
Each subsection ends with a new formula for the index
which the next  subsection builds on until we eventually
obtain, in subsection 5.5, part $2)$ of the main theorem. In 
Subsection 5.6 we will prove the index 
formula in part $3)$ of Theorem \ref{SFLIT}.
Our starting point is the McKean-Singer formula (Corollary 3.17) for the index
while in \cite{CPRS2} the starting point 
was the spectral flow formula of Carey-Phillips
\cite{CP2}.

\subsection{Exploiting Clifford-Bott periodicity}\label{getz}

We utilise
an idea of Getzler from \cite{G}
adapted to a more functional analytic setting
based on \cite{CP0}. 
We begin with an even semifinite spectral triple
$(\A,\HH,\D)$
with ${\Z}_2$-grading $\gamma$.
We will assume that this spectral triple is
$n-$summable for any $n>q$ with $q\geq 1$ fixed once and for all.
If $p\in \A$ then our aim is to derive from McKean-Singer
a new formula for the 
index of $p\D^+p=p_{-}\D^+p_+$ where $\D^+=(1-\gamma)\D(1+\gamma)/4=
P^{\perp}\D P$ and $p_+=PpP$ and $p_{-}=P^{\perp}pP^{\perp}.$
(Note that what follows differs significantly from what is done in \cite{CPRS2})

\begin{defn}
Form the Hilbert space
$\tilde{\mathcal H}=\IC^2\otimes\mathcal H$ on which acts the semifinite
von Neumann algebra $\tilde{\mathcal N}=M_2({\C})\otimes\mathcal N$.
Introduce the two dimensional Clifford algebra in the form
$$\sigma_1= \left(\begin{array}{cc}
                   0 & 1 \\
                   1 & 0
                   \end{array} \right),\ \  \sigma_2= \left(\begin{array}{cc}
                   0 & -i \\
                   i & 0
                   \end{array} \right), \ \ \sigma_3=
 \left(\begin{array}{cc} 1 & 0 \\
                         0 & -1
\end{array} \right).$$
Let $1_2$ denote the $2\times 2$ identity matrix and define
the grading in $\tilde{\mathcal N}$ by
$\tilde\gamma = \sigma_3\otimes \gamma$
and a Clifford element $\tilde{\sigma}_2=\sigma_2\otimes 1\in\tilde{\mathcal N}$
which anticommutes with $\tilde\gamma$
where $1$ is the identity operator in $\mathcal N$.
\end{defn}

Let $p\in \mathcal A$ be a projection. 
Introduce the following operators affiliated to $\tilde{\mathcal N}$
on $\tilde{\mathcal H}$:
$$\tD = \sigma_3\otimes \D,\ \ \ \D_p=p\D p+(1-p)\D(1-p)= \D+[\D,p](1-2p),$$
$$\D_w= (1-w)\D + w(p\D p+(1-p)\D(1-p))=(1-w)\D+w\D_p=\D+w[\D,p](1-2p),$$
and noting that $\tilde{\sigma}_2 \left(1_2\otimes(2p-1)\right)=
\sigma_2\otimes(2p-1)$,
we define:
$$\tD_{w,s} = \sigma_3\otimes \D_w + s(\sigma_2\otimes(2p-1))=:
 \tD_w+s(\sigma_2\otimes(2p-1)),
 \ \  
w\in [0,1], s\in (-\infty,\infty).$$
Note that  $\sigma_3\otimes \D_w$ is odd (i.e., anticommutes with the
grading $\tilde\gamma$)
and that $\tilde\sigma_2$ and $\sigma_3\otimes \D_w$ anticommute.
Notice that 
$\frac{d}{ds}\tD_{w,s}=\sigma_2\otimes(2p-1)$.

We extend the trace $\tau$ on $\mathcal N$ to $\tau_2:=Tr_2\otimes\tau$ on 
$\tilde{\mathcal N}$ by taking
the matrix trace $Tr_2$ in the first tensor factor. 
There is a  graded Clifford trace (super trace) 
on $\tilde{\mathcal N}$ 
which we write as $S\tau(T)=\frac{1}{2}\tau_2((\sigma_3\otimes 1)\tilde\gamma 
T)$,
$T\in \tilde{\mathcal N}$,
and note that this reduces to $\tau(\gamma S)$ for 
$T=1_2\otimes S\in\tilde{\mathcal N}.$

Now
\bean\tD_{w,s}^2&=& 1_2\otimes\D_w^2- s\sigma_3\sigma_2\otimes(2p-1)\D_w
+s\sigma_3\sigma_2\otimes\D_w(2p-1)+s^2\nno
&=& 1_2\otimes\D_w^2+ 2s(1-w)\sigma_3\sigma_2\otimes[\D,p]+s^2.\eean

Here we used $ \D_w(2p-1)-(2p-1)\D_w = 2(1-w)[\D,p].$
At $w=0$ we have
$$\tD^2_{0,s}=\left(\begin{array}{cc}
                   \D^2+s^2 & -i2s[\D,p] \\
                   -i2s[\D,p] & \D^2+s^2
                   \end{array} \right) =\tD^2+s^2+2s\s_3\s_2\otimes[\D,p]$$
and at $w=1$:
$$\tD^2_{1,s}=\left(\begin{array}{cc}
                   (p\D p+(1-p)\D(1-p))^2+s^2 & 0 \\
                   0 & (p\D p+(1-p)\D(1-p))^2+s^2 
                   \end{array} \right)=\tD_p^2+s^2, $$
where, $\tD_p:=\s_3\otimes\D_p.$
Note that
$$\frac{d}{dw}\D_w = p\D p+ (1-p)\D(1-p)-\D=[\D,p](1-2p).$$

\begin{lemma} \label{exact}
Consider the affine space $\Phi$ of perturbations, $\hat{\D}$, 
of $\tD=\s_3\otimes\D$ given by 
$$\Phi=\{\hat{\D}=\tD+X \ \ |\  X \in\tilde{\mathcal N}_{sa}\; and
\  [X,\sigma_2\otimes\gamma]=0\}$$
Notice that each $\hat{\D}$ commutes with $\sigma_2\otimes\gamma$.
Let, for any $n>q$ and $\hat{\D}\in\Phi$
$$\alpha_{\hat{\D}}(Y)=\tau_2\left((\sigma_2\otimes\gamma)
Y(1+\hat{\D}^2)^{-n/2}\right)$$
Then $\hat{\D}\mapsto \alpha_{\hat{\D}}$ is an exact one-form (i.e., an exact
section of the cotangent bundle to $\Phi$).
\end{lemma}

The proof of this Lemma is a trivial variation of the proof of 
Lemma 5.6 of \cite{CPRS2}: in the notation of that lemma, let 
$\Gamma=\s_2\otimes\gamma$  and multiply
$S\tau$ by $2$.

Now, for $n>q+1$ we introduce the function
$$ a(w) = \frac{1}{4}\int_{-\infty}^{\infty}\tau_2\left((1_2\otimes
\{\gamma (2p-1)\})(1+\tD_{w,s}^2)^{-n/2}\right)ds$$
This integral converges absolutely due to the following two 
estimates. The first is from Lemma 5.2 of \cite{CPRS2} (together with the
Remark immediately preceding that Lemma) with $n=q+2r$ and all $s\geq 2$:
$$\n(1+\tD_{w,s}^2)^{-n/2}\n_1=\n(1+\tD_{w,s}^2)^{-q/2-r}\n_1
\leq\n(1/2+\tD_w^2)^{-(q/2+\epsilon)}\n_1\left(
1/2+(1/2)(s^2)\right)^{-r+\epsilon}$$
where $r=\frac{n-q}{2}>\frac{1}{2}+\epsilon.$

The second is from Corollary 8 of Appendix B of \cite{CP1}, letting
$\tD_w=\tD+wA\in \Phi$ where $A$ is in $\tilde{\mathcal N}_{sa}.$ The cited 
result gives us a constant $C=C(A,q,\epsilon)>0$ such that
$$\n(1/2+\tD_w^2)^{-(q/2+\epsilon)}\n_1\leq
\tau_2\left((2(1+\tD_w^2)^{-1})^{(q/2+\epsilon)}\right)
\leq C\tau_2\left((1+\tD^2)^{-(q/2+\epsilon)}\right).$$
So, with $n=q+2r$ and $r>1/2$, the function $a(w)$ is well-defined.

With that settled, we now observe that 
\bean a(w)&=&\frac{1}{4}\int_{-\infty}^{\infty}\tau_2\left((\sigma_2\otimes
\gamma)(\s_2\otimes (2p-1))
(1+\tD_{w,s}^2)^{-n/2}\right)ds\nno
&=&
\frac{1}{4}\int_{-\infty}^{\infty}\tau_2\left((\sigma_2\otimes\gamma)
\left(\frac{d}{ds}\tD_{w,s}\right)(1+\tD_{w,s}^2)^{-n/2}\right)ds\eean
with the last expression designed to link with the result of the previous 
lemma.
In fact, $a(w)$ does not really depend on $w$ as we now prove.

\begin{lemma} We have that $a(w)$ is constant, in particular, $a(0)=a(1).$

\end{lemma}

\begin{proof}
Exactness of the one-form $\alpha$ means that integral of $\alpha$ along any
continuous piecewise smooth closed path in $\Phi$ must be $0$. Consider 
the closed (rectangular) path $\beta$ given by the four linear paths 
beginning with:
\bean\beta_{0,N}(s)=\tD_{0,s}\;\; for\;\; s\in [-N,N]; &\text {then}&
\beta_N(w)=\tD_{w,N}\;\; for\;\; w\in [0,1];\;\;\; \text {then}\nno
\beta_{1,N}(s)=\tD_{1,-s}\;\; for\;\; s\in [-N,N]; &\text{then}&
\beta_{-N}(w)=\tD_{1-w,-N}\;\; for \;\;w\in [0,1].\eean
Then, the integral of $\alpha$ around $\beta$ is $0$. For example,
the integral of
$\alpha$ along $\beta_N$ is:
$$\int_0^1\tau_2\left((\s_2\otimes\gamma)\frac{d}{dw}(\tD_{w,N})(1+\tD_{w,N}^2)
^{-n/2}\right)dw=\int_0^1\tau_2\left(B(1+\tD_{w,N}^2)^{-n/2}\right)dw,$$
where $B=\s_2\s_3\otimes\gamma[\D,p](1-2p)\in \tilde{\mathcal N}$. Now by the 
above estimates we have:
$$\n B(1+\tD_{w,N}^2)^{-n/2}\n_1
\leq C||B||\tau_2\left((1+\tD^2)^{-(q/2+\epsilon)}\right)
\left(1/2+(1/2)(N^2)\right)^{-r+\epsilon},$$
which for $r>\epsilon$ goes to $0$ as $N\to\infty.$ Similarly, the integral
along $\beta_{-N}$ goes to $0$ as $N\to\infty.$\\
Now the integral of $\alpha$ along $\beta_{0,N}$ is:
$$\int_{-N}^{N}\tau_2\left((1_2\otimes
\gamma (2p-1))(1+\tD_{0,s}^2)^{-n/2}\right)ds \to 4a(0)\;\; as\;\;
N\to\infty.$$
Similarly, the integral of $\alpha$ along $\beta_{1,N}$ converges to
$-4a(1)$ as $N\to\infty$. That is, $4a(0)-4a(1)=0$ or $a(0)=a(1).$
Similarly, $a(0)=a(w)$ for any $w\in [0,1].$
\end{proof}

Using the preceding lemma we obtain
$$a(1)=a(0)=\frac{1}{4}\int_{-\infty}^{\infty}
\tau_2\left((1_2\otimes\gamma(2p-1))(1+\tD_{0,s}^2)^{-n/2}\right)ds$$
and thus we can calculate $a(0)$ and $a(1)$ to obtain two different
expressions for the same quantity.

For the next calculation, observe that the definition of $a(w)$ gives
\bean a(1)=
\frac{1}{4}\int_{-\infty}^\infty\tau_2\left((1_2\otimes
\gamma(2p-1))(1+\tD_{1,s}^2)^{-n/2}\right)ds\eean
and inserting $\tD_{1,s}^2=\tD^2_p+s^2$,  we get by an application of 
McKean-Singer (Corollary 3.17):
\bean a(1)&=& \frac{1}{4}\int_{-\infty}^{\infty}
\tau_2\left((1_2\otimes\gamma(2p-1))(1+(p\tD p+(1-p)\tD(1-p))^2+s^2)^{-n/2}
\right)ds\nno
&=& \int_{-\infty}^{\infty}\tau
\left(\gamma p(1+(p\D p)^2+s^2)^{-n/2}\right)ds\nno
&&-\frac{1}{2}\int_{-\infty}^{\infty}
\tau\left(\gamma(1+(p\D p+(1-p)\D(1-p))^2+s^2)^{-n/2}\right)ds\nno
&=&\,Ind(p\D^+p)\int_{-\infty}^{\infty}(1+s^2)^{-n/2}ds-\frac{1}{2}
\int_{-\infty}^{\infty}
\tau\left(\gamma(1+\D_p^2+s^2)^{-n/2}\right)ds\eean

We put Lemma \ref{exact} to work again to get rid of the subscript $p$
in the last integral above.

\begin{lemma}
With the hypotheses as above and $n=q+2r>q+1$, we have:
$$\int_{-\infty}^{\infty}
\tau\left(\gamma(1+\D_p^2+s^2)^{-n/2}\right)ds=\int_{-\infty}^{\infty}
\tau\left(\gamma(1+\D^2+s^2)^{-n/2}\right)ds.$$
\end{lemma}

\begin{proof}
For $w\in [0,1]$ and $s\in \R$ we let:
$$\hat{\D}_{w,s} = \tD_w + s(\s_2\otimes 1)
                          = \tD + wA + s(\s_2\otimes 1)$$
where $A=\s_3\otimes([\D,p](1-2p)),$ so that $\tD+A=\tD_p.$ 
Then both perturbations of $\tD$ commute 
with $\s_2\otimes\gamma$ and therefore $\hat{\D}_{w,s}\in\Phi$.
Moreover $\s_2\otimes 1$ anticommutes with $\tD$ and with $A$, and so
$\hat{\D}_{w,s}^2 = \tD_w^2 + s^2.$
In particular, 
$\hat{\D}_{0,s}^2 = \tD^2 + s^2$ and
$\hat{\D}_{1,s}^2 = \tD_p^2 + s^2.$
One now applies Lemma \ref{exact} to the closed rectangular path in $\Phi$
described as follows: 
\bean\beta_{0,N}(s)=\hat{\D}_{0,s}\;\; for\;\; s\in [-N,N]; &\text {then}&
\beta_N(w)=\hat{\D}_{w,N}\;\; for\;\; w\in [0,1];\;\;\; \text {then}\nno
\beta_{1,N}(s)=\hat{\D}_{1,-s}\;\; for\;\; s\in [-N,N]; &\text{then}&
\beta_{-N}(w)=\hat{\D}_{1-w,-N}\;\; for \;\;w\in [0,1].\eean 
As in the previous lemma, the integral of the one-form along
$\beta_N$ equals:
$$\int_0^1\tau_2\left((\s_2\otimes\gamma)A(1+\hat{\D}_{w,N}^2)^{-n/2}\right)dw$$
and converges to $0$ as $N\to\infty.$ Similarly, the integral
along $\beta_{-N}$ goes to $0$ as $N\to\infty.$ 

Moreover, the integral along $\beta_{0,N}$ equals:
$$\int_{-N}^{N}\tau_2\left((\s_2\otimes\gamma)(\s_2\otimes 1)
(1+\tD^2+s^2)^{-n/2}\right)ds=2\int_{-N}^{N}\tau\left(\gamma
(1+\D^2+s^2)^{-n/2}\right)ds,$$ 
which as $N\to\infty$ converges to
$2\int_{-\infty}^{\infty}\tau\left(\gamma
(1+\D^2+s^2)^{-n/2}\right)ds.$
Similarly, the integral along $\beta_{1,N}$ converges to:
$-2\int_{-\infty}^{\infty}\tau\left(\gamma
(1+\D_p^2+s^2)^{-n/2}\right)ds.$
The proof is completed by observing that the integral around the closed 
path $\beta$ is $0$.
\end{proof}

This establishes the main formula of this section:
\begin{lemma}
For $n=q+2r>q+1$ we have:
$$
Ind(p\D^+p)C_{n/2}\;=\;a(1)+\frac{1}{2}\int_{-\infty}^{\infty}\tau
\left(\gamma(1+\D^2+s^2)^{-n/2}\right)ds$$
\bean 
&=&\frac{1}{4}\int_{-\infty}^{\infty}
\tau_2\left((1_2\otimes\gamma(2p-1))(1+\tD_{0,s}^2)^{-n/2}\right)ds+
\int_{0}^{\infty}\tau
\left(\gamma(1+\D^2+s^2)^{-n/2}\right)ds
\label{aofone}\eean
where
\bean C_{n/2}=\int_{-\infty}^\infty(1+s^2)^{-n/2}ds=\frac
{\Gamma(1/2)\Gamma(n/2-1/2)}{\Gamma(n/2)}.\eean
\end{lemma}
 
Note that $C_{n/2}$ is the normalisation `constant' that appeared in 
\cite{CPRS2}. Given the expression in terms of $\Gamma$ functions
we may take $n$ as a complex variable 
and see that the first pole is at $n=1$. 
If we write $n=q+2r$ then the pole is at $r=(1-q)/2$
which is the origin of the critical point
in the zeta functions in our main theorem. We reiterate that the
the above formula is only valid for $r>1/2$ but that the LHS gives an analytic
continuation of the RHS to a deleted neighbourhood of this critical point
$r=(1-q)/2$.

\subsection{Resolvent Expansion of the Index}\label{resolventexp}

In this subsection we will take the index formula of the preceding lemma
and apply a resolvent expansion to the integrand.
We begin with some notation.
Let $N=[(q+1)/2]$, where $[\cdot]$ denotes the integer part. If $q$ is an even 
integer, then $N=q/2$. If $q$ is an odd integer, then $N=(q+1)/2$. In general,
since $N\leq (q+1)/2 <N+1$ we have $2N-1\leq q <2N+1,$ so that $2N-1$ is the 
greatest odd integer in $q$. Also, $2N\leq q$ whenever $2n\leq q<2n+1$ for some 
positive integer $n$. In all cases $2N+2>q$.

We allow $q\geq 1$, so $(1+\D^2)^{-n/2}\in\LL^1(\mathcal N)$ for all $n>q$. 
By scale invariance of the index, we may replace $\D$ by $\epsilon\D$ 
without changing the index. Since we need $\n[\D,p]\n<\sqrt{2}$ below, we 
assume this without further comment.
We now make use of the Clifford structure. It allows us to employ the 
resolvent expansion to study $a(0)$, and we need only retain the even terms.

\begin{lemma}
Let $l$ be the line $\{\lambda=a+iv:-\infty<v<\infty\}$ where $0<a<1/2$ is 
fixed. There exists  $1>\delta>0$ such that for $r>1/2$
\bea &&a(1)=a(0)=\frac{1}{2}\int_{-\infty}^\infty 
S\tau(\{1_2\otimes(2p-1)\}(1+\tD_{0,s}^2)^{-q/2-r})ds\nno
&=&\frac{1}{2}\int_{-\infty}^\infty S\tau\left(\frac{1}{2\pi i}\int_l
\lambda^{-q/2-r}(1_2\otimes(2p-1))
(\lambda-(1+\tD^2_{0,s}))^{-1}d\lambda\right)ds
\nno
&=&\frac{1}{2}\int_{-\infty}^\infty S\tau\left(\frac{1}{2\pi i}\int_l
\lambda^{-q/2-r}(1_2\otimes(2p-1))\sum_{m=0}^{2N}\left(
(\lambda-(1+s^2+\tD^2))^{-1}2s\s_3\s_2\otimes[\D,p]\right)^m\times\right.\nno
&&\qquad\qquad\qquad\qquad\left.(\lambda-(1+s^2+\tD^2))^{-1}
d\lambda\right)ds+holo\eea
where $holo$ is a function of $r$ holomorphic for $Re(r)>(1-q)/2-\delta/2$.
\end{lemma}

\begin{proof} The first equality is just Cauchy's formula
$f(z)=\frac{1}{2\pi i}\int_l f(\lambda)(\lambda-z)^{-1}d\lambda$
(see the introductory remarks of section 6.2 of \cite{CPRS2} addressing
the issue of convergence).
The expansion in the statement of the Lemma is just the resolvent expansion:
$$\tilde{R}_s(\lambda)=\sum_{m=0}^{2N}\left(R_s(\lambda)2s\s_3\s_2\otimes
[\D,p]\right)^mR_s(\lambda) + \left(R_s(\lambda)2s\s_3\s_2\otimes[\D,p]
\right)^{2N+1}\tilde{R}_s(\lambda),$$
where $\tilde{R}_s(\lambda)=
(\lambda-(1+s^2+\tD^2+2s\s_3\s_2\otimes[\D,p]))^{-1}$ and $R_s(\lambda)=
(\lambda-(1+s^2+\tD^2))^{-1}$.
The remainder term in the resolvent expansion is 
\be \frac{1}{2}\int_{-\infty}^\infty (2s)^{2N+1}S\tau\left(\frac{1}{2\pi i}
\int_l\lambda^{-q/2-r}(1_2\otimes(2p-1))\left(\s_3\s_2\otimes 
R_s(\lambda)[\D,p]
\right)^{2N+1}\tilde{R}_s(\lambda)d\lambda\right)ds.
\label{rem}\ee
By H\"{o}lder's inequality
\ben \n (R_s(\lambda)[\D,p])^{2N+1}\n_1\leq \n[\D,p]\n^{2N+1}
\n R_s(\lambda)\n_{2N+1}^{2N+1}=C\n 
R_s(\lambda)^{2N+1}\n_1,\een
and by \cite[Lemma 5.3]{CPRS2} for all sufficiently small $\epsilon>0$
and $q\geq1$.
\ben \n R_s(\lambda)^{2N+1}\n_1\leq C_\epsilon
((1/2+s^2-a)^2+v^2)^{-(2N+1)/2+(q+\epsilon)/4}\een
where $\lambda=a+iv$. Moreover for $\n[\D,p]\n<\sqrt{2}$ we have by 
\cite[Lemma 5.1]{CPRS2}
\ben \n \tilde{R}_s(\lambda)\n\leq 
C'((1+s^2-a-s\n[\D,p]\n)^2+v^2)^{-1/2}.\een
We put these estimates together to obtain an estimate 
for the trace norm of the remainder term (\ref{rem}). We find
\bean \n(\ref{rem})\n_1&\leq &C''_\epsilon\int_{-\infty}^\infty s^{2N+1}
\int_{-\infty}^\infty \sqrt{a^2+v^2}^{-q/2-r}
((1/2+s^2-a)^2+v^2)^{-(2N+1)/2+(q+\epsilon)/4}\times
\nno
&\times&((1+s^2-a-s\n[\D,p]\n)^2+v^2)^{-1/2}dvds.\eean
Applying \cite[Lemma 5.4]{CPRS2} (one easily checks that one can
integrate from $-\infty$ instead of $0$ there)
we find that this integral is finite provided 
$q+\epsilon<2N+2$ and $-2N-2r+\epsilon<0$. The first condition is always 
satisfied by virtue of our choice of $2N$ and $\epsilon\leq 1$. 
For the second condition to be true 
at $q+2r=1-\delta$ requires that $\epsilon+\delta+q<2N+1$
and a $\delta$ satisfying this condition can always be found since
$2N-1\leq q < 2N+1$. That (\ref{rem})
defines a holomorphic function of $r$ for $Re(r)>(1-q)/2-\delta/2$ can
be seen by an argument essentially identical to the one in the proof of
\cite[Lemma 7.4]{CPRS2}. 
\end{proof}

{\bf Observation} Since $1_2\otimes\gamma$ commutes with $\tD^2$ and
anticommutes with $1_2\otimes[\D,p]$, all the terms in the expansion 
with $m$ odd vanish. On the other hand, each of the integrands with $m$
even is an even function of $s$ and so we may replace 
$\frac{1}{2}\int_{-\infty}^\infty$ by $\int_0^\infty$ in the above expansion. 

{\bf Observation} Using \cite[Lemma 7.2]{CPRS2}, we find that for $Re(r)>0$
each term in the 
above sum is in fact trace class, so we may interchange the trace and the sum. 
Having done this, we examine the $m=0$ term.
The $m=0$ term in the above expansion is given by 
\ben 2 \int_0^\infty\tau(\gamma p(1+s^2+\D^2)^{-q/2-r})ds-\int_{0}^{\infty}
\tau\left(\gamma
(1+\D^2+s^2)^{-q/2-r}\right)ds,\een
where the second term is the same (except for sign) as the second
term in Lemma 6.5. 
Hence if we write $R_s(\lambda)=(\lambda-(1+s^2+\tD^2))^{-1}$ 
we have for $r>1/2$
\bea &&Ind(p\D^+p)C_{q/2+r}=2\int_0^\infty\tau(\gamma p(1+s^2+\D^2)^{-q/2-r})ds
\;\;+\label{resexp}\\
&&\sum_{m=2,even}^{2N}\int_0^\infty\hspace{-.1in} S\tau
\left(\frac{1}{2\pi i}\int_l
\lambda^{-q/2-r}(1_2\otimes(2p-1))
\left(R_s(\lambda)2s\s_3\s_2\otimes[\D,p]\right)^mR_s(\lambda)d\lambda\right)ds+
holo\nonumber\eea
where $holo$ is a function of $r$ holomorphic for $Re(r)>(1-q)/2-\delta/2$.

The left hand side of Equation (\ref{resexp}) provides an
analytic continuation of the right hand side which is otherwise only defined
for $Re(r)>1/2$. 
The simple pole at $r=(1-q)/2$ has residue equal to $Ind(p\D^+p)$. 
We intend to compute this residue in terms of the analytic continuations
of the integrals appearing on 
the right hand side.

\subsection{Pseudodifferential Expansion of the Index}\label{psidoexp}

In this section we use ideas of \cite{CPRS2} and Connes-Moscovici's
pseudodifferential calculus to rewrite equation (\ref{resexp})
 in a form in which
all the resolvents in the integrand are commuted to the right.
In this new form we will be in a position to calculate residues
explicitly term by term. Our aim is to prove the following:

\begin{lemma} There exists a $1>\delta>0$ such that for $r>1/2$  
$$Ind(p\D^+p)C_{q/2+r}
=2\int_0^\infty\tau(\gamma p(1+s^2+\D^2)^{-q/2-r})ds\hspace{.1in}
+\sum_{m=2,even}^{2N}\sum_{|k|=0}^{2N-m}(-1)^{m/2}C(k)\hspace{.05in}\times$$
$$\times\int_0^\infty\hspace{-.1in} (2s)^mS\tau\left(\frac{1}{2\pi i}
\int_l\lambda^{-q/2-r}\left(1_2\otimes\{(2p-1)[\D,p]^{(k_1)}\cdots
[\D,p]^{(k_m)}\}\right)R_s(\lambda)^
{|k|+m+1}d\lambda\right)ds+holo$$
where $R_s(\lambda)=(\lambda-(1+s^2+\tD^2))^{-1}$, $holo$ is a function of 
$r$ holomorphic for $Re(r)>(1-q)/2-\delta/2$ and
$C(k)={(|k|+m)!}\alpha(k)$.
\end{lemma}

\begin{proof}
This is an application of our adaptation of Higson's version of the 
pseudodifferential expansion, 
and the observation that $(\s_3\s_2)^2=-1$. 
By \cite[Lemma 6.11]{CPRS2}, the remainder from the pseudodifferential 
expansion (applied to the $m$-th term in the resolvent expansion) is of order 
at most $-2m-(2N-m)-3=-m-2N-3$. By  \cite[Lemma 6.12]{CPRS2}, the remainder 
$P_{m,N}$ satisfies
\ben \n P_{m,N}(s,\lambda)(\lambda-(1+s^2+\tD^2))^{(m+2N+3)/2}\n\leq C\een
where the bound is uniform in $s,\lambda$ and square roots use the principal 
branch of $\log$. We use this to replace 
$P_{m,N}$ by powers of the 
resolvent to estimate the trace norm of the remainder. We obtain
\bean &&\n\int_0^\infty s^m\int_l\lambda^{-q/2-r}P_{m,N}d\lambda ds\n_1\nno
&\leq&C\int_0^\infty s^m\int_{-\infty}^\infty\sqrt{a^2+v^2}^{-q/2-r}\n 
R_s(\lambda)^{(m+2N+3)/2}\n_1dvds\nno
&\leq&C'\int_0^\infty s^m\sqrt{a^2+v^2}^{-q/2-r}\sqrt{(1/2+s^2-a)^2+v^2}^
{-(m+2N+3)/2+(q+\epsilon)/2}dvds\eean
where the final estimate comes from \cite[Lemma 5.3]{CPRS2}. 
Applying \cite[Lemma 5.4]{CPRS2} we find that this integral is finite 
provided $m-2((m+2N+3)/2-(q+\epsilon)/2)<-1$ and 
$m-2((m+2N+3)/2-(q+\epsilon)/2)+1-q-2r<-2$. The former condition requires 
$q+\epsilon<2+2N$, which is true 
by our choice of $N$. For the second condition to be true at $q+2r=1-\delta$ 
requires $\epsilon+\delta+q<2N+1$, and since $2N-1\leq q <2N+1$, for 
sufficiently small $\epsilon>0$ there 
exists a $1>\delta>0$ satisfying this condition.
\end{proof}

\subsection{Integrating Out the Parameter Dependence}\label{int}

The formula of the last lemma has
two integrals: one over the resolvent parameter $\lambda$ 
and the other over $s\in[0,\infty)$. 
The $\lambda$ integral can be performed by 
a simple application of Cauchy's formula for derivatives.

\begin{lemma} There exists $1>\delta>0$ such that for $r>1/2$
\bean &&Ind(p\D^+p)C_{q/2+r}-2\int_0^\infty\tau(\gamma 
p(1+s^2+\D^2)^{-q/2-r})ds\nno
&=&\sum_{m=2,even}^{2N}\sum_{|k|=0}^{2N-m}(-1)^{m/2+|k|}C(k)
\frac{\Gamma(q/2+r+|k|+m)}{\Gamma(q/2+r)(|k|+m)!}\times\nno
&\times&\int_0^\infty\hspace{-.1in}(2s)^mS\tau\left(1_2\otimes
\{(2p-1)[\D,p]^{(k_1)}
\cdots[\D,p]^{(k_m)}\}(1+s^2+\tD^2)^{-q/2-r-|k|-m}\right)ds+holo\eean
where $holo$ is a function of $r$ holomorphic for $Re(r)>(1-q)/2-\delta/2$.
\end{lemma}

\begin{proof} After ``pulling'' the unbounded operator 
$1_2\otimes\{(2p-1)[\D,p]^{(k_1)}\cdots[\D,p]^{(k_m)}\}$ out of the integral
(how to do this is explained in the proof of \cite[Lemma 7.2]{CPRS2}) we just
apply Cauchy's Formula in the operator setting (also discussed in
\cite[Lemma 7.2]{CPRS2}):
\bean \frac{1}{2\pi i}\int_l\lambda^{-q/2-r}R_s(\lambda)^{|k|+m+1}
d\lambda&=&\frac{1}{(|k|+m)!}\left.\left(\frac{d^{|k|+m}}{d\lambda^{|k|+m}}
\lambda^{-q/2-r}\right)\right|_{\lambda=(1+s^2+\tD^2)}\eean
\hspace*{2.9in}$=(-1)^{|k|+m}\frac{\Gamma(q/2+r+|k|+m)}
{\Gamma(q/2+r)(|k|+m)!}(1+s^2+\tD^2)^{-q/2-r-|k|-m}.$
\end{proof}

The remaining $s$-integral is not difficult either. 

\begin{lemma}\label{integrated} There exists $1>\delta>0$ such that for $r>1/2$
\bean &&Ind(p\D^+p)C_{q/2+r}-C_{q/2+r}\tau(\gamma p(1+\D^2)^{(1-q)/2-r})\nno
&=&\sum_{m=2,even}^{2N}\sum_{|k|=0}^{2N-m}(-1)^{m/2+|k|}C(k)
\frac{2^{m-1}\Gamma((m+1)/2)
\Gamma(q/2+r+|k|+(m-1)/2)}{\Gamma(q/2+r)(|k|+m)!}\times\nno
&&\qquad\qquad\times S\tau\left(1_2\otimes\{(2p-1)[\D,p]^{(k_1)}\cdots
[\D,p]^{(k_m)}\}
(1+\tD^2)^{-q/2-r-|k|-(m-1)/2}\right)\nno
&=&\sum_{m=2,even}^{2N}\sum_{|k|=0}^{2N-m}(-1)^{m/2+|k|}C(k)
\frac{2^{m-1}\Gamma((m+1)/2)
\Gamma(q/2+r+|k|+(m-1)/2)}{\Gamma(q/2+r)(|k|+m)!}\times\nno
&&\qquad\qquad\times \tau\left(\gamma(2p-1)[\D,p]^{(k_1)}\cdots[\D,p]^{(k_m)}
(1+\D^2)^{-q/2-r-|k|-
(m-1)/2}\right)+holo\eean
where $holo$ is a function of $r$ holomorphic for $Re(r)>(1-q)/2-\delta/2$.
\end{lemma}

\begin{proof} The integral is a Bochner integral (for a discussion of the 
subtleties see the proof of \cite[Proposition 8.2]{CPRS2}), and so  
we can move the $s$-integral past the supertrace. Then using the Laplace
Transform argument of \cite[Proposition 8.2]{CPRS2}, we have:
\bean &&\int_0^\infty (2s)^m(1+s^2+\tD^2)^{-|k|-m-q/2-r}ds\nno
&=&\frac{1}{\Gamma(q/2+r+m+|k|)}\int_0^\infty\int_0^\infty 
u^{|k|+m+q/2+r-1}(2s)^me^{-(1+\tD^2)u}e^{-s^2u}dsdu\nno
&=&\frac{\Gamma((m+1)/2)2^{m}}{2\Gamma(|k|+m+q/2+r)}\int_0^\infty u^{|k|+
(m-1)/2+q/2+r-1}e^{-(1+\tD^2)u}du\nno
&=&\frac{2^{m-1}\Gamma((m+1)/2)\Gamma(|k|+(m-1)/2+q/2+r)}{\Gamma(|k|+m+q/2+r)}
(1+\tD^2)^{-|k|-(m-1)/2-q/2-r}.\eean
Substituting this into the result of the last lemma almost gives the result. The
only extra things to do are to note the value of the constant arising from the 
integration for $m=0$ and to  trace out the Clifford variables (which could 
have been done earlier). Removing the Clifford variables is easy, because it 
is just a trace over the $2\times 2$ identity matrix, and the factor of $1/2$ 
in the definition of the super trace cancels it out. Hence the result.
\end{proof}

\subsection{Simplifying the Constants}\label{consts}

To obtain the constants that appear before the residues of the zeta functions
in the statement of our main theorem
requires us to  manipulate the constants in front of the zeta functions
in the statement of the last lemma of the preceding subsection.
Legendre's duplication formula for the Gamma function \cite[p. 200]{A}
says
$$2^{m-1}\Gamma((m+1)/2)=\sqrt{\pi}{\Gamma(m)}/{\Gamma(m/2)}.$$
For $m=0$ replace the right hand side with $\sqrt{\pi}/2$. Since for $m> 0$
and even, 
$\frac{\Gamma(m)}{\Gamma(m/2)}=\frac{1}{2}\frac{m!}{(m/2)!},$
we have
$ 2^m\Gamma((m+1)/2))=\sqrt{\pi}{m!}/{(m/2)!}.$
The functional equation for the Gamma function says
\bean \frac{\sqrt{\pi}\Gamma(r+(q-1)/2+|k|+m/2)}{\Gamma(q/2+r)}&=&
\frac{\sqrt{\pi}\Gamma(r+(q-1)/2)}{\Gamma(q/2+r)}\prod_{j=0}^{|k|+m/2-1}
(r+(q-1)/2+j)\nno
&=& C_{q/2+r}\sum_{j=1}^{|k|+m/2}(r+(q-1)/2)^j\s_{(|k|+m/2),j},\eean
where the $\s_{(|k|+m/2),j}$ are the elementary symmetric functions
of the integers $1,2,...,|k|+m/2$. Substituting 
these oddments into the formula from Lemma \ref{integrated} 
for $r>1/2$ and with $h=m/2+|k|$  gives, {\em modulo 
functions of $r$ holomorphic for $Re(r)>(1-q)/2-\delta$}: 
 \bea &&Ind(p\D^+p)C_{q/2+r}-C_{q/2+r}\tau(\gamma p(1+\D^2)^{(1-q)/2-r})\nno
&=&\sum_{m=2,even}^{2N}\sum_{|k|=0}^{2N-m}(-1)^{m/2+|k|}C(k)
\frac{\sqrt{\pi}m!
\Gamma(q/2+r+|k|+(m-1)/2)}{2(m/2)!\Gamma(q/2+r)(|k|+m)!}\times\nno
&&\qquad\qquad\times \tau\left(\gamma(2p-1)[\D,p]^{(k_1)}\cdots[\D,p]^{(k_m)}
(1+\D^2)^{-q/2-r-|k|-(m-1)/2}\right)\nno
&=&\sum_{m=2,even}^{2N}\sum_{|k|=0}^{2N-m}(-1)^{m/2+|k|}
\frac{m!}{(m/2)!}\frac{\alpha(k)}{2}
\frac{\sqrt{\pi}\Gamma(q/2+r+|k|+(m-1)/2)}{\Gamma(q/2+r)}\times\nno
&&\qquad\qquad\times \tau\left(\gamma(2p-1)[\D,p]^{(k_1)}\cdots[\D,p]^{(k_m)}
(1+\D^2)^{-q/2-r-|k|-
(m-1)/2}\right)\nno
&=&\sum_{m=2,even}^{2N}\sum_{|k|=0}^{2N-m}(-1)^{m/2+|k|}
\frac{m!}{(m/2)!}\frac{\alpha(k)}{2}
C_{q/2+r}\times\nno
&&\times\sum_{j=1}^{h}\s_{h,j}(r+(q-1)/2)^j \tau\left(\gamma(2p-1)
[\D,p]^{(k_1)}\cdots[\D,p]^{(k_m)}
(1+\D^2)^{(1-q)/2-r-|k|-m/2}\right). \label{penultimate}\eea

Observe that we have not used the isolated spectral dimension
assumption at any point in this calculation. 
Despite this, the above sum of zeta functions (which includes the $m=0$
term which we have written once on the LHS of the first equality to save space) 
has a simple pole at $r=(1-q)/2$ with residue equal to $Ind(p\D^+p)$. 
This proves part $2)$ of Theorem \ref{SFLIT}.\\
\hspace*{.2in}To proceed further, we  need to assume that the individual 
zeta functions have analytic continuations.

\subsection{Taking the Residues}

This  step will prove the index formula in part $3)$ of Theorem \ref{SFLIT}.
We now have to {\bf assume isolated spectral dimension}. Then, denoting: 
$$\zeta_{m,k}(z)=\tau\left(\gamma(2p-1)[\D,p]^{(k_1)}\cdots[\D,p]^{(k_m)}
(1+\D^2)^{-z-|k|-m/2}\right),$$
we have for $r>1/2$
$$Ind(p\D^+p)C_{q/2+r}=C_{q/2+r}\tau(\gamma p(1+\D^2)^{-(q-1)/2-r})$$
$$+\hspace{-.1in}\sum_{m=2,even}^{2N}\sum_{|k|=0}^{2N-m}(-1)^{m/2+|k|}
\frac{m!}{(m/2)!}\frac{\alpha(k)}{2}
C_{q/2+r}\sum_{j=1}^{h}\s_{h,j}(r+(q-1)/2)^j\zeta_{m,k}((q-1)/2+r)
+holo$$
where $h=|k|+m/2$. 
Now, divide through by $C_{q/2+r}$, and multiply by $1/(r+(q-1)/2)$. 
The remainder term is now
\ben \frac{holo}{C_{q/2+r}(r+(q-1)/2)},\een
which is still holomorphic at the critical point 
(since it has a removable singularity).  
Denote the analytic continuation of $\zeta_{m,k}((q-1)/2+r)$ by 
${\bf Z}_{m,k}((q-1)/2+r)$. Define for $j=-1,0,1,...$
\ben \tau_j\left(\gamma(2p-1)[\D,p]^{(k_1)}\cdots
[\D,p]^{(k_m)}(1+\D^2)^{-m/2-|k|}
\right)=res_{(1-q)/2}(r-(1-q)/2)^j{\mathbf Z}_{m,k}((q-1)/2+r),\een
(we replace $2p-1$ by $2p$  when $m=0$). Thus taking the residues of the 
left and right hand sides of Equation 
(\ref{penultimate}) we obtain (setting $h=|k|+m/2$)
\bea &&Ind(p\D^+p)=res_{r=(1-q)/2}\frac{1}{(r+(q-1)/2)}Ind(p\D^+p)\nno
&=&res_{r=(1-q)/2}\frac{1}{(r+(q-1)/2)}\tau(\gamma p(1+\D^2)^{-(r-(1-q)/2)})\nno
&+&\sum_{m=2,even}^{2N}\sum_{|k|=0}^{2N-m}(-1)^{h}\frac{m!}{(m/2)!}\frac{
\alpha(k)}{2}\sum_{j=1}^{h}\s_{h,j}
res_{r=(1-q)/2}(r+(q-1)/2)^{j-1}{\bf Z}_{m,k}((q-1)/2+r)\nno
&=& \tau_{-1}(\gamma p)\eea
$$+\sum_{m=2,even}^{2N}\sum_{|k|=0}^{2N-m}(-1)^{h}\frac{m!}{(m/2)!}
\frac{\alpha(k)}{2}
\sum_{j=1}^{h}\s_{h,j}\tau_{j-1}\left(\gamma(2p-1)[\D,p]^{(k_1)}\cdots
[\D,p]^{(k_m)}(1+\D^2)^{-m/2-|k|}\right).\label{Indexpairing}$$

Observe that since $j-1$ runs from $0$ to $|k|+m/2-1$, at worst, we only need 
to consider the first $|k|+m/2-1$ terms in 
the principal part of the Laurent series for ${\bf Z}_{m,k}$ at $r=(1-q)/2$, as 
well as the constant term. Moreover, this number is bounded by
\ben |k|+m/2-1\leq 2N-m+m/2-1=2N-m/2-1\leq 2N-1\leq q\een
since $2N-1\leq q < 2N+1$. Hence $|k|+m/2-1\leq q$. Furthermore, since 

$\gamma(2p-1)[\D,p]^{(k_1)}\cdots[\D,p]^{(k_m)}\in OP^{|k|}$, it equals 
$B(1+\D^2)^{|k|/2}$ for some $B$ bounded, and so:
\ben \gamma(2p-1)[\D,p]^{(k_1)}\cdots[\D,p]^{(k_m)}
(1+\D^2)^{-m/2-|k|-r-(q-1)/2}=B(1+\D^2)^{-m/2-|k|/2-r-(q-1)/2}\een
The right hand side has finite trace for 
$$Re(r)>({1-m-|k|})/{2}=(1-q)/{2}+(q-m-|k|)/{2}.$$
Thus whenever $m+|k|>q$ we obtain a term which is holomorphic at $r=(1-q)/2$. If
$[q]$ is {\bf odd} then there exists
$n\in{\N}$ with $2n+1\leq q< 2n+2$ and so $N=n+1$ and $2N=2n+2>q$. Hence the 
residues of the terms with $m=2N$ all vanish, 
and similarly for any $m=2,...,2N-2$ the residues of the top terms 
with $|k|=2N-m$ vanish.\\
\hspace*{.2in}This computation, which has produced 
Equation (\ref{Indexpairing}) has actually 
proved the index formula in part $3)$ of Theorem \ref{SFLIT}.
To prove $1)$ and the remainder of $3)$, we need to study the resolvent cocycle.

\section{The Resolvent Cocycle in the Even Case}\label{resolventcocycle} 

Part $3)$ of  Theorem \ref{SFLIT} claims that
the index is actually a pairing of
a $(b,B)$ cocycle with the Chern character of the idempotent $p$. Similarly, in 
$1)$ we have an `almost' cocycle, and the residue of the pairing computes the 
index.
In order to show
this
we introduce an auxiliary function-valued $(b,B)$-cochain called the
{\bf resolvent cocycle} (cf \cite[Section 7]{CPRS2}). 
The definition is inspired by the resolvent expansion, 
and we show that it is a $(b,B)$-cocycle {\em modulo functions of $r$ 
holomorphic in an open half-plane containing $r=(1-q)/2$}. 
We use the resolvent cocycle to complete the proof of Theorem \ref{SFLIT} in 
subsection \ref{finalproof}.\\
\hspace*{.2in}Our starting point for this section is the expansion of $a(0)$ 
obtained in 
equation (\ref{resexp}) at the end of Subsection \ref{resolventexp}. We have
\bean &&Ind(p\D^+p)C_{q/2+r}=2\int_0^\infty\tau(\gamma p(1+s^2+\D^2)^{-q/2-r})
ds\nno
&+&\sum_{m=2,even}^{2N}\int_0^\infty S\tau\left(\frac{1}{2\pi i}\int_l
\lambda^{-q/2-r}(2p-1)\left(R_s(\lambda)2s\s_3\s_2\otimes[\D,p]
\right)^mR_s(\lambda)d\lambda\right)ds+holo\eean
where $holo$ is a function of $r$ holomorphic for $r>(1-q)/2-\delta/2$ where 
$1>\delta>0$.
If we now perform the `super bit' of the trace we obtain
\bean &&Ind(p\D^+p)C_{q/2+r}=2\int_0^\infty\tau(\gamma p(1+s^2+\D^2)^{-q/2-r})
ds\nno
&+&\sum_{m=2,even}^{2N}(-1)^{m/2}\int_0^\infty 
\tau\left(\frac{1}{2\pi i}\int_l
\lambda^{-q/2-r}\gamma(2p-1)\left(R_s(\lambda)2s[\D,p]\right)^mR_s(\lambda)
d\lambda\right)ds+holo\eean
where by abuse of notation we have written $R_s(\lambda)=
(\lambda-(1+s^2+\D^2))^{-1}$ (as opposed to $\tD^2$).

Assuming that the right hand side is (almost) the pairing of  a 
cocycle with the Chern character of the projection $p$, 
to obtain a formula for the cocycle we expect to 
remove the normalisations coming from the Chern character of 
$p$, and that is all. Including the powers of two in the normalisation gives 
the next definition.

\begin{defn} For $m$ even, $0\leq m\leq 2N$, $a_0,...,a_m\in\A$, and 
$\eta_m=2^{m+1}\frac{(m/2)!}{m!}$ define the following function of $r$ for 
$r>(1-m)/2$:
\ben \phi^r_m(a_0,...,a_m)=\frac{\eta_m}{2\pi i}
\int_0^\infty s^m
\tau\left(\gamma\int_l\lambda^{-q/2-r}a_0R_s(\lambda)
[\D,a_1]R_s(\lambda)\cdots R_s(\lambda)[\D,a_m]R_s(\lambda)d\lambda\right)ds.
\een
\end{defn}

Observe that the definition for $m=0$ and the Cauchy formula gives
\ben \phi^r_0(a_0)=2\int_0^\infty\tau(\gamma a_0(1+s^2+\D^2)^{-q/2-r})ds.\een

\begin{prop} For $0\leq m\leq 2N-2$, $B\phi^r_{m+2}+b\phi^r_m=0$ and there 
exists a $1>\delta>0$ such that $(b\phi^r_{2N})(a_0,...,a_{2N+1})$ is 
holomorphic for $Re(r)>(1-q)/2-\delta/2$.
\end{prop}

\begin{proof} We use the discussion of \cite[Subsection 7.2]{CPRS2}, noting 
some minor differences which arise due to the grading, $\gamma$. 
We begin by computing 
$B\phi^r_{m+2}$. Applying the definitions we have, 
\ben (B\phi^r_{m+2})(a_0,...,a_{m+1})
=\sum_{j=0}^{m+1}(-1)^j\phi^r_{m+2}(1,a_j,...,a_{m+1},a_0,...,a_{j-1})=\een
\ben\sum_{j=0}^{m+1}
\frac{(-1)^j\eta_{m+2}}{2\pi i}
\int_0^\infty\hspace{-.15in} s^{m+2}\tau\left(\gamma\int_l
\lambda^{-q/2-r}R_s(\lambda)[\D,a_j]
\cdots [\D,a_{m+1}]R_s(\lambda)\cdots [\D,a_{j-1}]R_s(\lambda)d\lambda\right)ds
\een
\ben=\sum_{j=0}^{m+1}\frac{\eta_{m+2}}{2\pi i}
\int_0^\infty s^{m+2}\tau\left(\gamma\int_l\lambda^{-q/2-r}[\D,a_0]
R_s(\lambda)\cdots R_s(\lambda)1R_s(\lambda)\cdots [\D,a_{m+1}]R_s(\lambda)
d\lambda\right)ds\een
The last line follows from \cite[Lemma 7.7]{CPRS2} modified by the fact that 
while $R_s(\lambda)$ commutes with $\gamma$, $[\D,a_i]$ 
anticommutes with $\gamma$.
We now employ \cite[Lemma 7.6]{CPRS2}:
\ben -k\int_0^\infty s^{k-1}\frac{1}{2\pi i}\tau\left(\gamma\int_l
\lambda^{-q/2-r}A_0R_s(\lambda)A_1R_s(\lambda)\cdots A_mR_s(\lambda)
d\lambda\right)ds\een
\ben=2\sum_{j=0}^{m}\int_0^\infty s^{k+1}\frac{1}{2\pi i}\tau\left(
\gamma\int_l\lambda^{-q/2-r}A_0R_s(\lambda)A_1R_s(\lambda)\cdots A_j
R_s(\lambda)1R_s(\lambda)A_{j+1}\cdots A_mR_s(\lambda)d\lambda\right)ds.\een
Applying this formula to our computation for $B\phi^r_{m+2}$ yields
\bean &&(B\phi^r_{m+2})(a_0,...,a_{m+1})\nno
&=&-\frac{1}{2}(m+1)\frac{\eta_{m+2}}{2\pi i}
\int_0^\infty s^m\tau\left(\gamma\int_l\lambda^{-q/2-r}R_s(\lambda)[\D,a_0]
\cdots R_s(\lambda)[\D,a_{m+1}]R_s(\lambda)d\lambda\right)ds\nno
&=&-\frac{\eta_m}{2\pi i}\int_0^\infty s^m\tau\left(\gamma\int_l
\lambda^{-q/2-r}R_s(\lambda)[\D,a_0]\cdots[\D,a_{m+1}]
R_s(\lambda)d\lambda\right)ds.\eean
Here we used $(m+1)/2\times\eta_{m+2}=\eta_m$.
Next one expands the first commutator on the right hand side, $[\D,a_0]=
\D a_0-a_0\D$, and anticommutes the second $\D$ through the remaining
$[\D,a_j]$ using $\D[\D,a_j]+[\D,a_j]\D=[\D^2,a_j].$ Recalling that $\D$ 
anticommutes with $\gamma$, we find from the proof of 
\cite[Proposition 7.10]{CPRS2} that
\ben (B\phi^r_{m+2})(a_0,...,a_{m+1})=\een
\ben\frac{\eta_m}{2\pi i}\int_0^\infty\hspace{-.1in}s^m
\sum_{j=1}^{m+1}(-1)^{j+1}
\tau\left(\gamma\hspace{-.05in}\int_l\lambda^{-q/2-r}
R_s(\lambda)a_0R_s(\lambda)[\D,a_1]\cdots[\D^2,a_j]\cdots [\D,a_{m+1}]
R_s(\lambda)d\lambda\right)ds.\een
We recall that $B\phi^r_0=0$, by definition.
The computation of $b\phi^r_m$ is precisely the same as 
\cite[Proposition 7.10]{CPRS2}, and gives
\ben (b\phi^r_{m})(a_0,...,a_{m+1})=\een
\ben\frac{\eta_m}{2\pi i}\int_0^\infty s^m\sum_{j=1}^{m+1}(-1)^j\tau\left(
\gamma\int_l\lambda^{-q/2-r}
R_s(\lambda)a_0R_s(\lambda)[\D,a_1]\cdots[\D^2,a_j]\cdots [\D,a_{m+1}]
R_s(\lambda)d\lambda\right)ds.\een
Hence $B\phi_{m+2}^r+b\phi^r_m=0$ for $0\leq m\leq 2N-2$ (indeed for all $m\geq 
0$).

For $m=2N$, we use H\"{o}lder's inequality 
(together with \cite[Lemma 6.10]{CPRS2} 
to see that $|[\D^2,a_j]R_s(\lambda)|\leq C^{\prime} |R_s(\lambda)|^{1/2}$)
which yields a constant $C$ independent of $s$ and $\lambda$ so that:
\ben \n R_s(\lambda)a_0R_s(\lambda)[\D,a_1]\cdots 
R_s(\lambda)[\D^2,a_j]R_s(\lambda)\cdots 
R_s(\lambda)[\D,a_{2N+1}]R_s(\lambda)\n_1\leq C\n R_s(\lambda)^{2N+5/2}\n_1.
\een

Consequently, we have the estimate (using \cite[Lemma 5.3]{CPRS2})
$$|(b\phi^r_{2N})(a_0,...,a_{2N+1})|
\leq C\int_0^\infty s^{2N}\int_{-\infty}^\infty\sqrt{a^2+v^2}^{-q/2-r}\n 
R_s(\lambda)^{2N+5/2}\n_1dvds$$
\bean 
&\leq & C_\epsilon\int_0^\infty s^{2N}
\int_{-\infty}^\infty\sqrt{a^2+v^2}^{-q/2-r}
\sqrt{(1/2+s^2-a)^2+v^2}^{-2N-5/2+(q+\epsilon)/2}dvds.\eean
Consulting  \cite[Lemma 5.4]{CPRS2} we find that this integral is convergent 
when $2r=1-q-\delta$ provided $q+\epsilon-4<2N$, which is true, and 
$q+\epsilon+\delta<2N+3$, which again is true. As for the case of the
remainder term in the proof of 
\cite[Lemma 7.4]{CPRS2} this shows that the above formula for 
$(b\phi^r_{2N})(a_0,...,a_{2N+1})$
gives a holomorphic function of $r$ in a neighbourhood of $(1-q)/2$ as claimed.
\end{proof}
Observe that together with Equation (\ref{resexp}) the above result proves 
part $1)$ of Theorem \ref{SFLIT}.

\subsection{The Residue Cocycle}\label{finalproof}
In this subsection we complete the proof of Theorem \ref{SFLIT}.
First we need to define the residue cocycle.
\begin{defn} Let $(\A,\HH,\D)$ be a $QC^\infty$ finitely summable spectral 
triple with isolated spectral dimension $q\geq 1$. For $m=2,...,2N$ and 
$a_0,...,a_m\in\A$ define functionals
\bean&&\phi_m(a_0,...,a_m)\nno
&=&\sum_{|k|=0}^{2N-m}(-1)^{|k|}
\alpha(k)\sum_{j=1}^{|k|+m/2}\s_{(|k|+m/2),j}\tau_{j-1}
\left(\gamma a_0[\D,a_1]^{(k_1)}\cdots[\D,a_m]^{(k_m)}(1+\D^2)^{-m/2-|k|}\right)
,\eean
and for $m=0$ define $\phi_0(a_0)=\tau_{-1}(\gamma a_0).$
\end{defn}

\begin{thm} Let $(\A,\HH,\D)$ be a $QC^\infty$ finitely summable spectral 
triple with isolated spectral dimension $q\geq 1$. When evaluated on any 
$a_0,...,a_m\in\A$, the components 
$\phi^r_m$ of the resolvent cocycle $(\phi^r)$ analytically continue to 
a deleted neighbourhood of
$r=(1-q)/2$.
 Moreover, if we denote this continuation by $\mathbf{\Phi}_m^r(a_0,...,a_m)$ 
 then
\ben res_{r=(1-q)/2}\frac{1}{C_{q/2+r}(r+(q-1)/2)}\mathbf{\Phi}^r_m(a_0,...,a_m)
=\phi_m(a_0,...,a_m).\een
\end{thm}

{\bf Remark} Observe that, as a function of $r$,
$[C_{q/2+r}(r+(q-1)/2)]^{-1}$
has a removable singularity at $r=(1-q)/2$. 
Thus all the statements concerning the resolvent cocycle also apply 
to the resolvent cocycle multiplied by this function.

\begin{proof}  
For $m$ even, evaluate $\phi^r_m$ on $a_0,...,a_m\in\A$ and apply the 
pseudodifferential 
expansion. 
This yields (modulo functions holomorphic for $Re(r)>(1-q)/2-\delta$)
$$\phi^r_m(a_0,...,a_m)=$$
$$\sum_{|k|=0}^{2N-m}C(k)\frac{\eta_m}{2\pi i}
\int_0^\infty s^m\tau\left(\gamma\int_l\lambda^{-q/2-r}a_0[\D,a_1]^{(k_1)}
\cdots[\D,a_m]^{(k_m)}R_s(\lambda)^{m+|k|+1}d\lambda\right)ds.$$
Proceeding according to our previous computations we have
$$\phi^r_m(a_0,...,a_m)$$
\bean &=&\sum_{|k|=0}^{2N-m}\frac{(-1)^{m+|k|}C(k)\Gamma(q/2+r+m+|k|)}
{\Gamma(q/2+r)(m+|k|)!}2^{m+1}\frac{(m/2)!}{m!}\times\nno
&\times&\int_0^\infty s^m\tau\left(\gamma a_0[\D,a_1]^{(k_1)}\cdots
[\D,a_m]^{(k_m)}(1+s^2+\D^2)^{-(m+|k|+q/2+r)}
\right)ds\ \ \ \ \ \mbox{Cauchy}\nno
&=&\sum_{|k|=0}^{2N-m}\frac{(-1)^{|k|}C(k)2^m\Gamma((m+1)/2)
\Gamma(q/2+r+m/2-1/2+|k|)}{\Gamma(q/2+r)(m+|k|)!}\frac{(m/2)!}{m!}\times\nno
&\times&\tau\left(\gamma a_0[\D,a_1]^{(k_1)}\cdots[\D,a_m]^{(k_m)}
(1+\D^2)^{-(m/2-1/2+|k|+q/2+r)}\right)
\ \ \ \ \qquad\qquad s-\mbox{integral}\nno
&=&\sum_{|k|=0}^{2N-m}\frac{(-1)^{|k|}C(k)\sqrt{\pi}
\Gamma(q/2+r+m/2-1/2+|k|)}{\Gamma(q/2+r)(m+|k|)!}\times\nno
&\times&\tau\left(\gamma a_0[\D,a_1]^{(k_1)}\cdots
[\D,a_m]^{(k_m)}(1+\D^2)^{-(m/2-1/2+|k|+q/2+r)}
\right)\ \ \ \qquad\qquad m-\mbox{factorials}\nno
&=&\sum_{|k|=0}^{2N-m}(-1)^{|k|}C_{q/2+r}\alpha(k)\sum_{j=0}^{h}\s_{h,j}
(r+(q-1)/2)^j\times\nno
&\times&\tau\left(\gamma a_0[\D,a_1]^{(k_1)}\cdots[\D,a_m]^{(k_m)}
(1+\D^2)^{-(m/2-1/2+|k|+q/2+r)}
\right)\ \ \ \Gamma\ \mbox{function and constants},
\eean
where $h=|k|+m/2$.
The result is now clear.
\end{proof}

\begin{cor} Let $(\A,\HH,\D)$ be a $QC^\infty$ finitely summable spectral 
triple with isolated spectral dimension $q\geq 1$. The cochain $(\phi)$ with 
components $\phi_m$, $m=0,2,...,2N$, is a $(b,B)$-cocycle. For any projection 
$p\in\A$ we have
\ben Ind(p\D^+p)=\sum_{m=0,even}^{2N}\phi_m(Ch_m(p)).\een
\end{cor}

\begin{proof} The first statement follows because $\left(
(C_{q/2+r}(r+(q-1)/2))^{-1}\phi^r\right)$ is a $(b,B)$-cocycle modulo 
functions holomorphic at $r=(1-q)/2$ and hence so is its analytic 
continuation,\\
$\left((C_{q/2+r}(r+(q-1)/2))^{-1}\mathbf{\Phi}^r\right)$. 
For the second statement we recall that
\ben Ch_0(p)=p,\ \ \ Ch_m(p)=(-1)^{m/2}\frac{m!}{2(m/2)!}(2p-1)\otimes 
p^{\otimes m}.\een
Thus $\sum_{m=0}^{2N}\phi_m(Ch_m(p))$ is given precisely by the formula on the 
right hand side of Equation (\ref{Indexpairing}), the left hand side of which 
is $Ind(p\D^+p)$. This completes the proof.
\end{proof}

We have now completed the proof of Theorem \ref{SFLIT}. We present the easy 
proof of Corollary \ref{lowdim}.

\begin{cor} For $1\leq q <2$, we do not need to assume
isolated spectral dimension to compute the index pairing. 
\end{cor}

\begin{proof} For $1\leq q<2$ we have $N=1$, but as we observed after Equation 
(\ref{Indexpairing}), the term with $m=2N$ is holomorphic at $r=(1-q)/2$ when 
$[q]$ is an odd integer. Hence we have only the $m=0$ term. So,
\ben Ind(p\D^+p)=\tau(\gamma p(1+\D^2)^{-(q-1)/2-r})+\frac{holo}{C_{q/2+r}}.\een
By the Remark in Theorem 6.4:
\ben \frac{1}{r+(q-1)/2}Ind(p\D^+p)=\frac{1}{r+(q-1)/2}
\tau(\gamma p(1+\D^2)^{-(q-1)/2-r})+holo.\een
Taking residues we have
\ben Ind(p\D^+p)=res_{r=(1-q)/2}\tau(\gamma p(1+\D^2)^{-(q-1)/2-r}),\een
and the residue on the right necessarily exists, and is equal to 
$\tau_{-1}(\gamma p)$. Hence the individual terms in the expansion of the 
index analytically continue to a punctured neighbourhood of $r=(1-q)/2$ with 
no need to invoke the isolated spectral dimension hypothesis. The single 
term $\phi_0$ forms a $(b,B)$ cocycle for $\A$ since $B\phi_0=0$ and 
$b\phi^r_0$ is holomorphic at $r=(1-q)/2$. 
\end{proof}

\end{document}